\documentclass[11pt,a4paper,english]{article}

\usepackage{mathptmx}
\usepackage{a4wide}
\usepackage{amssymb,amsthm}
\usepackage{amsmath}
\usepackage{amsfonts}
\usepackage{latexsym}
\usepackage{amsxtra}
\usepackage{graphicx}
\usepackage{wrapfig}
\usepackage{mathrsfs}
\usepackage{epsfig}
\usepackage[pdftex=true,breaklinks=true,bookmarks=true,linktocpage=true,pdfpagelabels=true,plainpages=false,pagebackref=true]{hyperref}

\hypersetup{
pdfauthor={Strzelecki,von der Mosel},
pdftitle={Tangent-point self-avoidance energies for curves},
pdfsubject={article,paper},
pdfproducer={LaTeX with hyperref package},
pdfkeywords={Menger curvature, curves},
hypertexnames=true,
pdfstartview={FitV},
pdfstartpage={1},
pdfhighlight={/I},
pdfborder={0 0 0 0},
bookmarksopen=false,
bookmarksnumbered=true,
pdfpagemode={UseOutlines},
unicode=false,
breaklinks=true,
frenchlinks=false,
a4paper=true,
pdftoolbar=true,
pdfmenubar=true,
pdfcenterwindow=true,
pdffitwindow=true,
pdfnewwindow=true,
colorlinks=true,
linkcolor=blue,
citecolor=blue,
urlcolor=blue,
filecolor=blue
}

\newcommand{\heikodetail}[1]{}
\renewcommand{\proof}{\noindent{\sc Proof:}\hskip 1.1em}
\renewcommand{\qed}{\hfill\mbox{$\Box$}\\}
\def\yy{\mbox{$\spadesuit$}}
\newcommand{\invisible}[1]{\par ($\spadesuit$ \emph{hidden comments in \TeX{}
file\/})}

\newcommand{\omitted}[1]{}

\newtheorem{lemma}{\bf Lemma}[section]
\newtheorem{theorem}[lemma]{\bf Theorem}
\newtheorem{proposition}[lemma]{\bf Proposition}
\newtheorem{corollary}[lemma]{\bf Corollary}

\newtheorem{definition}[lemma]{\bf Definition}

\newcommand{\Fo}{\,\,\,\text{for }\,\,}
\newcommand{\Foa}{\,\,\,\text{for all }\,\,}

\newcommand{\AND}{\,\,\,\text{and }\,\,}

\newcommand\Reals{{\mathbb R}}
\newcommand\R{{\mathbb R}}
\newcommand\Z{{\mathbb Z}}

\newcommand\N{{\mathbb N}}
\renewcommand\S{{\mathbb S}}

\newcommand{\bbbr}{\Reals}

\newcommand\dist{\mathop{\rm dist}\nolimits}
\newcommand\diam{\mathop{\rm diam}\nolimits}

\newcommand\ang{\mathop{\mbox{$<\!\!\!)$}}\nolimits}

\newcommand{\xx}{\mbox{$\clubsuit$}}

\newcommand{\E}{\mathcal{E}}

\newcommand{\A}{\mathcal{A}}
\newcommand{\M}{\mathcal{M}}

\newcommand{\eps}{\epsilon}
\renewcommand{\H}{\mathscr{H}}
\def\osc{\mathop{\rm osc\,}}         % oscillation  %

\begin{document}

%%%%%%%%%%%%%%%%%%%%%%%%%%%%%%%%%%%%%%%%%%%%%%%%%%%%
%%%%%%%%%%%   title   %%%%%%%%%%%%%%%%%%%%%%%%%%%
%%%%%%%%%%%%%%%%%%%%%%%%%%%%%%%%%%%%%%%%%%%%%%%%
\title{Tangent-point self-avoidance energies for curves}

\author{Pawe\l{} Strzelecki\thanks{PS and his research is partially
supported by the Polish Ministry of Science and Higher Education
grant no. N~N201 397737 (years 2009-2012).},
\setcounter{footnote}{6}\quad Heiko von der
Mosel\thanks{HvdM is partially supported by the DFG grant
Mo966/4-1.}}

\maketitle
%%%%%%%%%%%%%%%%%%%%%%%%%%%%%%%%%%%%%%%%%%%%%%
%%%%%%%%%%%% abstract  %%%%%%%%%%%%%%%%%%%%%%%
%%%%%%%%%%%%%%%%%%%%%%%%%%%%%%%%%%%%%%%%%%%%%

\frenchspacing

\begin{abstract}

We study a two-point self-avoidance energy $\E_q$ which is defined
for all recti\-fiable curves in $\R^n$ as the double integral
along the curve of $1/r^q$. Here $r$ stands for the radius of the
(smallest) circle that is tangent to the curve at one point and
passes through another point on the curve, with obvious natural
modifications of this definition in the exceptional, non-generic
cases. It turns out that finiteness of $\E_q(\gamma)$ for $q\ge 2$
guarantees that $\gamma$ has no self-intersections or triple
junctions and therefore must be homeo\-morphic to the unit circle
$\S^1$ or to a closed interval $I$.  For $q>2$ the energy $\E_q$
evaluated on curves in $\R^3$ turns out to be a knot energy
separating different knot types by infinite energy barriers and
bounding the number of knot types below a given energy value. We
also establish an explicit upper bound on the Hausdorff-distance
of two curves in $\R^3$ with finite $\E_q$-energy that guarantees
that these curves are ambient isotopic. This bound depends only on
$q$ and the energy values of the curves. Moreover, for all $q$
that are larger than the critical exponent $q_{\text{crit}}=2$,
the arc\-length parametrization of $\gamma$\, is of class
$C^{1,1-2/q}$, with H\"{o}l\-der norm of the unit tangent
depending only on $q$, the length of $\gamma$, and the local
energy. The exponent $1-2/q$ is optimal.

\vspace{3mm}

\centering{Mathematics Subject Classification (2000): 28A75,
49Q10, 53A04, 57M25}

\end{abstract}

\renewcommand\theequation{{\thesection{}.\arabic{equation}}}
\def\setnumbers{\setcounter{equation}{0}}

%%%%%%%%%%%%%%%%%%%%%%%%%%%%%%%%%%%%%%%%%%%%%
%%%%%%%%%%%%  intro  = Section 1 %%%%%%%%%%%%%%%%%%%%%%%%%%
%%%%%%%%%%%%%%%%%%%%%%%%%%%%%%%%%%%%%%%%%%%%%%

\section{Introduction}\label{sec:1}
Imagine a space craft travelling with constant speed along an unknown
and possibly quite irregular closed path $\Gamma$
in an unexplored territory of the universe. After some time $L>0$ the
loop is completed at least once, and the only data the astronauts can
measure at time $t$ are the ratios of the squared distance from any
previous position $\Gamma(s)$, to the distance of the current line
of direction $\ell(t)$ from that previous position $\Gamma(s)$, i.e.,
the quotients
\begin{equation}\label{1.1}
2r(\Gamma(t),\Gamma(s)):=\frac{|\Gamma(t)-\Gamma(s)|^2}{\dist(\ell(t),\Gamma(s))}\in
[0,\infty]\quad\Fo s<t.
\end{equation}
What can the astronauts say about their path of travel? In other words,
how much information about a closed curve of finite length in Euclidean
space is encoded in the relative tangent-point data \eqref{1.1}?
The answer is: If the astronauts
obtain a finite integral mean of some
inverse power of all these data (after time $2L$) they can extract essential
topological information as well as explicit smoothness properties of their path
of travel!

To make this precise we assume from now on that the path $\Gamma\subset\R^n$ is
a rectifiable curve of finite length, parametrized by arclength on the
circle $S_L\cong\R/(L\Z)$ of perimeter $L$. Hence,  $\Gamma$ is a
(not necessarily injective)
Lipschitz continuous mapping with $|\Gamma'|=1$ a.e. on $S_L$.
Geometrically, the tangent-point function
$$
r(\Gamma(t),\Gamma(s))=\frac{|\Gamma(t)-\Gamma(s)|^2}{2\dist(\ell(t),\Gamma(s))}
=\frac{|\Gamma(t)-\Gamma(s)|}{2\sin\ang(\Gamma'(t),\Gamma(s)-\Gamma(t))},
$$
involving the tangent line $\ell(t):=\{\Gamma(t)+\mu\Gamma'(t):\mu\in\R\}$
and
defined for all $s\in S_L$ and almost all $t\in S_L$, determines the radius
of the unique circle that is
tangent to $\Gamma$ at the position $\Gamma(t)$ and
passes through $\Gamma(s)$. (This radius is set to be zero
if $\Gamma(t)=\Gamma(s)$, and is infinite
if the vector $\Gamma(s)-\Gamma(t)\not=0$ is parallel to the tangent $\Gamma'(t)$).

The only assumption in the result indicated above is finiteness of the
{\it tangent-point potential}
\begin{equation}\label{1.2}
\E_q(\Gamma):=\int_0^L\int_0^L\frac{dsdt}{r^q(\Gamma(t),\Gamma(s))}
\quad\textnormal{for some $q\ge 2.$}
\end{equation}
\begin{theorem}[Finite energy path is a manifold]\label{thm:1.1}
If $\E_q(\Gamma)<\infty$ for some $q\ge 2$ then the image $\Gamma(S_L)$
is a one-dimensional topological manifold (possibly with boundary),
embedded in $\R^n.$
\end{theorem}
In particular, the image curve has no self-intersections, although there is
no chance to deduce injectivity of the arclength parametrization $\Gamma$
itself, since the integrand depends only on the image $\Gamma(S_L).$
Take, for example, a $k$-times covered circle of length $L/k$,
for which the integrand is constant,
$r(\Gamma(t),\Gamma(s))\equiv r_0$ for all $s,t\in S_L$, so that the energy
amounts to
$$
\E_q(\textnormal{$k$-times covered circle})=\frac{L^2}{r_0^q}=
k^2\int_0^{L/k}\int_0^{L/k}\frac{dsdt}{r_0^q}=k^2\E_q(\textnormal{once-covered
circle})<\infty.
$$
So the space craft's course
cannot be too wild, since it traces a one-dimensional manifold without any
non-tangential self-crossings. But without further input the astronauts have
no clue of how often they have completed that course. Moreover, in case
their path forms
a manifold with boundary, say, a circular arc, there would be an abrupt
(and for the crew probably quite noticeable) change of direction at the endpoints
of that arc. Mathematically, one can easily reparametrize the manifold to obtain
a new injective arclength parametrization, which translates to the additional
information that the spacecraft  does not pass by any previous position at all,
$\Gamma(t)\not=\Gamma(s)$ for all $t\not=s,$ which we will assume from now on.

In light of Theorem \ref{thm:1.1} the tangent-point potential
$\E_q$ evaluated on closed curves in $\R^3$ may serve as a valid
{\it knot energy} as suggested by Gonzalez and Maddocks in
\cite[Section 6]{GM99}, that is, as a functional separating
different knot types by infinite energy barriers. It was shown by
Sullivan \cite[Prop. 2.2]{sullivan} that for $q>2$ the energy
$\E_q$ blows up on a sequence of smooth knots converging smoothly
to a smooth curve with self-crossings. (His proof uses the
Taylor formula up to order two for the converging curves, and a
uniform bound for the remainders.) As a consequence of our
analysis we generalize this result to continuous curves replacing
smooth convergence by uniform convergence (see Proposition
\ref{self-repulsive}). Thus $\E_q$ for $q>2$ is indeed {\it
self-repulsive} or {\it charge}, and hence a knot energy according
to the definition given by O'Hara \cite[Def. 1.1]{ohara-book},
which provides an affirmative answer to an open question posed in
\cite[Problem 8.1]{ohara-book}.
%In the scale
%invariant case $q=2$, as Gonzalez and Maddocks indicate, $\E_q$ is
%directly related to the knot energies considered by Buck and Simon
%\cite{bucksimon1997,bucksimon1999}.
It also turns out that $\E_q$ is \emph{strong\/} for $q>2$: among
all continuous closed curves $\gamma$ of fixed length $L$ and
$\E_q(\gamma)<E$ there are only finitely many knot types, see
Proposition~\ref{prop:strong}. This gives a partial answer to a
conjecture expressed by Sullivan in \cite[p. 184]{sullivan}
(leaving open the case $q=2$, and we do not consider links with
more than one component). Both these knot-theoretic results are
based on a priori $C^{1,\alpha}$-estimates for curves of finite
$\E_q$-energy, discussed in more detail later on.

We will show in addition that two curves, whose Hausdorff-distance
is bounded above by an explicit small constant depending only on
the energy values, are in fact in the same knot class. A
qualitative version of such an isotopy result is well-known in the
smooth category; see e.g. \cite[Chapter 8]{hirsch}, or
\cite{blatt}. Here, however, we have explicit quantitative bounds.
Notice also that Hausdorff-distance alone, no matter how small,
does not suffice to separate knot classes;\footnote{Consider for
example two different torus knots on the surface of a very thin
rotational torus; for the classification of torus knots see e.g.
Burde and Zieschang \cite[Chapter 3.E]{burde-zieschang}.} bounded
$\E_q$-energy is crucial here.
\begin{theorem}[Isotopy]\label{thm:1.2}
For any $q>2$ there is an explicit constant $\delta(q)>0$ depending only
on $q$ such that any two closed rectifiable
curves with injective arclength parametrizations
$\Gamma_1,$ $\Gamma_2$, with finite $\E_q$-energy, are ambient isotopic if
their Hausdorff-distance is less than
$$
\delta(q)\max\{\E_q(\Gamma_1),\E_q(\Gamma_2)\}^{-\frac{1}{q-2}}.
$$
\end{theorem}
Our proof of Theorem \ref{thm:1.2} follows closely the arguments of
Marta Szuma\'{n}ska who proved in her Ph.D. thesis a similar result
\cite[Chapter 5]{marta} for a related three-point potential, the {\it
integral Menger curvature}
\begin{equation}\label{1.3}
\M_p(\Gamma):=\int_0^L
\int_0^L
\int_0^L\frac{dsdtd\sigma}{R^p(\Gamma(s),\Gamma(t),\Gamma(\sigma))},\quad p>3,
\end{equation}
where $R(x,y,z)$ denotes the circumcircle radius of three points
$x,y,z$ in Euclidean space. Essentially one reduces the isotopy
question to that between polygons inscribed in $\Gamma_1$ and
$\Gamma_2$, whose edge lengths are solely controlled in terms of
the energy. For polygonal knots a similar result is contained in
the work of Millet, Piatek, and Rawdon \cite[Theorem 4.2]{MPR},
where instead of \eqref{1.3}, the polygonal  thickness of the
polygons together with their edge length determines the smallness
condition on the Hausdorff distance that guarantees isotopy of two
polygonal knots. For general curves, thickness was defined by
Gonzalez and Maddocks in \cite{GM99} as the smallest possible
circumcircle  radius $R(\cdot,\cdot,\cdot)$ when evaluated on
all triples of distinct curve points. This concept of thickness
was used as a tool in variational applications involving curves
and elastic rods subject to various topological constraints; see
e.g. \cite{GMSvdM}, \cite{CKS}, \cite{heiko1}--\cite{heiko3},
\cite{GevdM09}, \cite{GevdM10}, and has been studied numerically,
\cite{CPR}, \cite{CLMS05}, \cite{ashton_etal}.

The inverse of thickness of a curve $\Gamma$
can be obtained as limits $\M_p^{1/p}(\Gamma)$ for $p\to\infty$,
or $\E_q^{1/q}(\Gamma)$ for $q\to\infty.$
In our papers \cite{svdm-single,ssvdm-double,ssvdm-triple} we have
studied regular\-izing, self-avoidance and compactness effects of
several integral energies, including $\M_p$, which involve, vaguely speaking, various
bounds for $1/R$ understood as a function of three variables,
including bounds in $L^p$, in $L^p(X_1,L^\infty(X_2))$ where
$X_1=S_L$ and $X_2=S_L\times S_L$ (or vice versa), and in spaces
that resemble the classic Morrey spaces $L^{p,\lambda}$. In each
case we were able to detect similar phenomena: there is a certain
limiting exponent for which an appropriate functional is scale
invariant, and above this exponent three sorts of effects take
place. First, curves with finite energy have no self-intersections.
Second, these energies serve well as knot energies allowing for valuable
compactness results for equibounded families of loops in fixed isotopy classes,
which is due to the third, the regularizing effect: Curves   with finite
energy are more regular than initially assumed.

For the present tangent-point potential $\E_q$ we obtain the following
regularity theorem, which shows that the astronauts would not experience
any sudden change of direction during their travel.
\begin{theorem}[Regularity]\label{thm:1.3}
If $q>2$ and the arclength parametrization
$\Gamma:S_L\to\R^n$ is chosen to be injective, then $\Gamma$
is continuously differentiable with a H\"older continuous tangent,
i.e., $\Gamma$ is of class  $C^{1,1-(2/q)}.$
More precisely, for each $q>2$ there exist two constants
$\delta(q)>0$ and $c(q)<\infty$ depending only on $q$  such that each injective
arclength parametrization $\Gamma$ with $\E_q(\Gamma)<\infty$
satisfies
\begin{equation}
|\Gamma'(u)-\Gamma'(v)| \le c(q) \left( \int_u^v\int_u^v
\frac{ds\, dt}{r(\Gamma(s),\Gamma(t))^q}\right)^{1/q}
|u-v|^{1-2/q} \label{ineq:1.2}
\end{equation}
for all $u,v\in S_L$ with $|u-v|\le \min \bigl(\delta(q)
\E_q(\Gamma)^{-1/(q-2)}, \frac 12 \diam \gamma\bigr)$.
\end{theorem}
The exponent $q=2$ is a limiting one here. It is relatively easy
to use scaling arguments and check that $\E_q(\Gamma)=\infty$ for
each $q\ge 2$ when $\Gamma$ parametrizes a closed polygonal curve,
but polygons have finite energy for all $q<2$. The resulting
H\"{o}l\-der exponent $1 - 2/q$ is reminiscent of the
classic Sobolev imbedding theorem in the supercritical case: the
domain of integration is two-dimensional, and the integrand is
related to curvature. For $C^2$-curves the behaviour of $1/r$
close to the diagonal of $S_L\times S_L$ (where $1/r$ might blow
up for curves with low regularity) encodes some information about
curvature, i.e. about second derivatives of the arclength
parametrization $\Gamma$. The point is that we need no information
about the existence of $\Gamma''$ in order to prove
 Theorem~\ref{thm:1.3}. A priori, we deal with curves that are
rectifiable only, and even the existence of $\Gamma'$ at all
parameters cannot be taken for granted.

Note that inequality \eqref{ineq:1.2} is qualitatively optimal:
for curves of class $C^{1,1}$ the integrand $1/r$ is bounded, and
\eqref{ineq:1.2} yields then $|\Gamma'(u)-\Gamma'(v)|\lesssim
|u-v|\sup (1/r) $ for $u,v$ sufficiently close; nothing stronger
can be expected as the familiar example of a stadium curve shows.
We discuss other examples briefly at the end of
Section~\ref{last}.

Before describing the main ideas of the proof and the structure
of the paper we would like
to mention
that while
working on generalizations of self-avoidance energies to surfaces in $\bbbr^3$,
see \cite{mengersurf}, which involved a search for suitable
inte\-grands, we have realized that $\E_q$ is a model energy that
might be the easiest one to extend to the fully general case, i.e.
to submanifolds of arbitrary dimension and co-dimension
\cite{tanpoint}. This was one of the motivations to write the
present note: to lay out in a simple, relatively easily tractable
case all the arguments that should be applicable in much greater
generality.

\medskip

Theorem~\ref{thm:1.1} is obtained as a corollary of
a slightly more general result, see Theorem~\ref{beta-mfd} below.
We first prove a technical lemma (see Section~\ref{beta-section})
which shows how $\E_q$ can be used to control the behaviour of the
so-called P. Jones' \emph{$\beta$-numbers\/},
\begin{equation}
\beta_\gamma(x,r) := \inf\biggl\{ \sup_{y\in \gamma \cap B(x,r)}
\frac{\dist (y,G)}{r} \quad \colon \quad \mbox{$G$\, is a straight
line through $x$} \biggr\}\, ,
\end{equation}
for small $r>0$ and closed balls $B(x,r)$ of radius $r$ with
center $x$. It turns out that if $\E_2(\Gamma)<\infty$ then
$\beta_\gamma(x,r)\to 0$ as $r\to 0$ uniformly with respect to
$x$, see Lemma~\ref{lem:2.3} and the remark at the end of
Section~\ref{beta-section}. And this is the key point to prove that
$\gamma=\Gamma(S_L)$ is a topological manifold, as we have the
following.
\begin{theorem}
\label{beta-mfd} If $\, \Gamma\colon S_L\to \bbbr^n$ is
arclength, and the image $\gamma=\Gamma(S_L) $ satisfies
\begin{equation}
\sup_{x\in\gamma} \beta_\gamma(x,d) \le \omega(d)\, \label{sup-b}
\end{equation}
where $\omega\colon [0,L]\to \bbbr$ is a continuous nondecreasing
function with $\omega(0)=0$, then $\gamma$ is a one-dimensio\-nal
sub\-manifold of $\bbbr^n$ (possibly with boundary).
\end{theorem}

The main idea behind the proof of Theorem~\ref{beta-mfd} is
simple: if the result were false, then we could find a point $x$
in $\gamma$ where a \emph{triple junction\/} occurs; in a small
ball $B$ centered at $x$  we would have (at least) three disjoint
arcs of $\gamma$ in a long narrow tube. Two of them would then be
very close (i.e., would leave $B$ crossing $\partial B$ in the
same spherical cap at one end of the tube). Observing points of
those two arcs, and using \eqref{sup-b} on smaller and smaller
scales, we are able to obtain a contradiction and eventually show
that there could be no triple junction at $x$. For details, see
Section~\ref{image}.

By the preliminary results of Section \ref{beta-section},
if $\E_q(\Gamma)<\infty$
for some $q>2$, then the control of $\beta$ numbers is much better
than just \eqref{sup-b}. Namely,
\begin{equation}
\sup_{x\in\gamma} \beta_\gamma(x,r)\lesssim r^\kappa \label{bk}
\end{equation}
for $\kappa=(q-2)/(q+4)< \lambda=1-2/q$; the constant in
\eqref{bk} depends on $\E_q(\Gamma)$. Applying \eqref{bk}
iteratively, we find in Section~\ref{diff} suitably defined cones
that contain short arcs of $\gamma$ and obtain an estimate for
their opening angles, proving that $\Gamma'$ exists everywhere and
is of class\footnote{Let us remark that for an $m$-dimensional set
$\Sigma\subset \bbbr^n$ that is \emph{Reifenberg flat with
vanishing constant\/} uniform estimates of $\beta$-numbers imply
that $\Sigma$ is a $C^{1,\kappa}$-manifold, see David, Kenig and
Toro \cite{dkt} and Preiss, Tolsa and Toro \cite{ptt}. Here, we
have no Reifenberg flatness a priori -- in general rectifiable
curves do not have to be Reifenberg flat; in fact we prove it by
hand, using energy bounds leading to \eqref{bk}.} $C^\kappa$.

Section \ref{sec:5} contains the proof of the isotopy result,
Theorem \ref{thm:1.2}. In the last section we show how to
bootstrap the initial gain of $C^{1,\kappa}$-regularity obtained
in Section~\ref{diff}, to the optimal regularity $\Gamma\in
C^{1,1-(2/q)}$, and we will establish \eqref{ineq:1.2}. We
stress the fact that Inequality \eqref{ineq:1.2} in
Theorem~\ref{thm:1.3} provides a uniform a priori estimate. This
can be used in variational applications and to ensure compactness
for infinite families of curves with uniformly bounded energy.
Some results of that type have been stated in
\cite{ssvdm-double,ssvdm-triple}; we do not follow that thread
here.

Finally, let us say that, at the moment, we have no clue how
$\Gamma'$ behaves in the limiting case $q=2$ (we do not even know
if it is defined everywhere for curves with finite $\E_2$-energy)
but we are tempted to think that $\Gamma'$ has vanishing mean
oscillation for $q=2$ and that local oscillations of the tangent
can be controlled by the local energy of the curve.

\subsection*{Notation}
We write $G(x,y)$ to denote the straight line through two
distinct points $x,y\in \R^n$. If $x=\Gamma(s),y=\Gamma(t)\in
\gamma:=\Gamma\bigl(S_L\bigr)\subset \bbbr^n$, then, abusing the
notation slightly, we write sometimes $G(s,t)$ instead of
$G(\Gamma(s), \Gamma(t))$.

For  a closed set $F$ in $\bbbr^n$ we set
\[
U_\delta(F) :=\{ x\in \bbbr^n\, \colon \, \dist(x,F)<\delta\},
\qquad \delta>0.
\]
In some places, it will be more convenient to work directly with
the slabs $U_\delta$ around appropriately selected lines than to
deal with the information expressed only in the language of
$\beta$-numbers. Finally, in Section~\ref{diff} we work with cones. For
$x\not=y\in\bbbr^n$ and $\eps\in (0,\frac\pi 2)$ we denote by
\begin{equation}
\label{cone-def} C_\eps(x;y)\, :=\, \{z\in\bbbr^n\colon \exists \
t\neq 0 \ \ \textrm{such that}  \  \ \ang(t(z-x),y-x) < \frac
\eps2 \}
\end{equation}
the double  cone whose  vertex is at the point $x$, with cone axis
passing through $y$, and with opening angle $\eps$.  All balls
$B(x,r)$ with radius $r>0$ and center $x\in\R^n$ are closed balls
throughout the paper.

\section{Decay of beta numbers}

\label{beta-section}

\begin{lemma}
\label{beta} Let $\E_q(\Gamma)$ be finite. There exists a constant
$c_0=c_0(q)>0$ such that if $\eps< 1/200$ and $d< \diam\gamma$
satisfy
\begin{equation}
\label{ed} \eps^{4+q}d^{2-q} \ge c_0(q) \E_q(\Gamma)\, ,
\end{equation}
then for every two points of the curve such that
$|\Gamma(s)-\Gamma(t)|= d$ we have
\[
\gamma\cap B_{2d} \bigl(\Gamma(s)\bigr) \ \subset \ U_{20\eps d}
\bigl( G(s,t)\bigr)\, .
\]
In particular,
\[
\beta_\gamma(\Gamma(s),2d) \le 10 \eps.
\]
\end{lemma}

For $q>2$ we set $\kappa= (q-2)/(q+4)$.

\begin{corollary}
There exists a $\delta_1=\delta_1(q)>0$ such that if
$\E_q(\Gamma)^{1/(q+4)} d^\kappa< \delta_1$, then
\[
\beta_\gamma(\Gamma(s), 2d) \le  c_1(q) \E_q(\Gamma)^{1/(q+4)}
d^\kappa\, .
\]
\end{corollary}

\omitted{By the results of David, Kenig and Toro \cite[Section
9]{davidkenigtoro} and Preiss, Tolsa and Toro \cite{ptt}, the
uniform decay of $\beta$ numbers for a curve $\gamma$ implies that
$\gamma$ is a $C^{1,\kappa}$-submanifold of $\bbbr^n$. Thus,
$\Gamma'\in C^\kappa$ once we know that $\Gamma$ is 1--1. (We
could also prove that directly, as in the triple energy paper;
iterations of the lemma give a result analogous to UCF.)}

\bigskip
\noindent\textbf{Proof of Lemma \ref{beta}.} For $s,t\in S_L $,
$d=|\Gamma(s)-\Gamma(t)|>0$ and $\eps>0$ small, we set
\begin{eqnarray*}
A_d(s,\eps) & := & \Gamma^{-1} (B_{\eps^2 d}(\Gamma(s)) \, = \, \{
\tau \in S_L\, \colon \, \Gamma(\tau) \in B_{\eps^2
d}(\Gamma(s))\}\, ,\\
X_d(s,t,\eps) & := & \{\sigma\in A_d(s,\eps)\, \colon \,
\Gamma'(\sigma) \,\textnormal{exists with}\,\, \ang
(\Gamma'(\sigma), \Gamma(t)-\Gamma(s)) \in \bigl[ \frac{\eps}{10},
\pi - \frac{\eps}{10}\bigr]\}\, , \\
N_d(s,t,\eps) & := & A_d(s,\eps) \setminus X_d(s,t,\eps)\, .
\end{eqnarray*}
Note that $|A_d(s,\eps)|\ge 2\eps^2 d$. The proof has two steps:
\begin{itemize}
\item we use the inequality
\[
\E_q(\Gamma) \ge \int_{X_d(s,t,\eps)}\int_{A_d(t,\eps)} r^{-q}\,
d\sigma\, d\tau
\]
to show that $X_d(s,t,\eps)$ must be a small subset of
$A_d(s,\eps)$, so that $|N_d(s,t,\eps)|\gtrsim \eps^2 d$;
\item we argue by contradiction, using energy estimates again, and
show the desired inclusion.
\end{itemize}

\medskip\noindent\textbf{Step 1.}
Fix $\sigma\in X_d(s,t,\eps)$ and $\tau \in A_d(t,\eps)$. We shall
show that $1/r((\Gamma(\sigma),\Gamma(\tau)) \gtrsim \eps/d$.

Since $|\Gamma(s)-\Gamma(t)|=d$, the triangle inequality yields
\begin{equation}
\label{stau} 2d> d(1+2\eps^2) \ge |\Gamma(\sigma)-\Gamma(\tau)|\ge
d (1-2\eps^2)\, .
\end{equation}
Let
\[
x:= \Gamma(\sigma) + d
\frac{\Gamma(t)-\Gamma(s)}{|\Gamma(t)-\Gamma(s)|}=\Gamma(\sigma) +
\bigl(\Gamma(t) -\Gamma(s) \bigr)\, \in \bbbr^n\, .
\]
Then $|x-\Gamma(\sigma)|=d$ and $|\Gamma(\tau)-x|\le 2\eps^2 d$ by
the triangle inequality. By definition of $X_d(s,t,\eps)$, the
angle $\alpha= \ang(x-\Gamma(\sigma),\Gamma'(\sigma))$ is
contained between $\eps/10$ and $\pi - \eps/10$. Therefore
\[
\dist(x,\ell(\sigma)) = d\sin \alpha \ge d \sin \frac{\eps}{10}
\ge \frac{d\eps}{20},
\]
and
\begin{equation}
\dist(\Gamma(\tau),\ell(\sigma))  \ge \frac{d\eps}{20} - 2\eps^2 d
\ge \frac{d\eps}{25} \label{tau-line}
\end{equation}
(here we use $\eps< 1/200$). Combining \eqref{stau} and
\eqref{tau-line}, we obtain
\[
\frac{1}{2r(\Gamma(\sigma),\Gamma(\tau))} \ge \frac{d\eps}{25}
(2d)^{-2} = \frac{\eps}{100 d}\, .
\]
Integration gives
\[
\E_q(\Gamma)\ge \int_{X_d(s,t,\eps)}\int_{A_d(t,\eps)} r^{-q}\,
d\tau\, d\sigma \ge \mathrm{const}\, \cdot\, |X_d(s,t,\eps)|
\eps^{2+q}d^{1-q},
\]
as $A_d(t,\eps) \ge 2 \eps^2d$. If $|X_d(s,t,\eps)| \ge \frac 12
\eps^2 d$, then
\[
\E_q(\Gamma)\ge \mathrm{const}(q)\, \cdot\, \eps^{4+q}d^{2-q},
\]
which gives a contradiction for an appropriate choice of $
c_0(q)>0$ in the lemma.

Thus, we have
\[
|X_d(s,t,\eps)| <  \frac 12 \eps^2 d \qquad\mbox{and}\qquad
|N_d(s,t,\eps)| > \frac 32 \eps^2 d.
\]

\medskip\noindent\textbf{Step 2.} Suppose now that $\Gamma(\tau)
\in B_{2d}(\Gamma(s)) \setminus U_{20\eps d}\bigl(G(s,t)\bigr)$.
Fix a $\sigma\in  N_d(s,t,\eps)$.

Since then $|\Gamma(\sigma)-\Gamma(s)|< \eps^2 d$ and the (acute)
angle between the vectors $\Gamma'(\sigma)$ and $\Gamma(t) -
\Gamma(s)$ is very close to $0$ or $\pi$ (the difference is at
most $\eps/10$), one can check that in fact
\[
\ell(\sigma) \cap B_{2d} (\Gamma(s)) \ \subset \ U_{\eps
d}(G(s,t)) \cap B_{2d} (\Gamma(s))\, .
\]
Therefore the distance from $\Gamma(\tau)$ to $\ell(\sigma)$ is at
least $19\eps d$. If $\tau_1 \in A_d(\tau,\eps)$, then
\[
\dist (\Gamma(\tau_1),\ell(\sigma)) \ge 19\eps d - \eps^2 d \ge 18
\eps d\, ,
\]
and
\[
\frac{1}{r(\Gamma(\sigma),\Gamma(\tau_1))} \ge \frac{18 \eps
d}{(3d)^2} > \frac{\eps}{d}\, .
\]
Integrating this inequality, we obtain
\[
\E_q (\Gamma) \ge \int_{N_d(s,t,\eps)}\int_{A_d(\tau,\eps)}
r^{-q}\, d\tau_1\, d\sigma\ge \frac{3}2 \eps^2 d \cdot 2 \eps^2 d
\cdot \left(\frac \eps d\right)^q  = 3 \eps^{4+q}d^{2-q}\, .
\]
Again, for an appropriate choice of $c_0(q)$ this gives a
contradiction with \eqref{ed}. \hfill $\Box$

\medskip

Since the assumption $q>2$ was not used at all in the proof of the
lemma, it is easy to check that the same reasoning that was used
to obtain \eqref{ed} gives in fact the following

\begin{lemma}\label{lem:2.3}
Assume that $q=2$ and $\E_2(\Gamma)<\infty$. Then there exists a
constant $c>0$ such that
\begin{equation}
\sup_{x\in\gamma} \beta_\gamma(x,d) \le c \omega_E(d), \qquad d\le
\diam \gamma,
\end{equation}
where
\begin{equation}
\omega_E(d) := \sup\biggl( \int_A \int_B \frac{ds\,
dt}{r\bigl(\Gamma(s),\Gamma(t)\bigr)^2} \biggr)^{1/6},
\end{equation}
the supremum being taken over all pairs of subsets $A,B\subset
S_L$ with $\H^1(A), \H^1(B) \le \frac{1}{100} d$.
\end{lemma}

\medskip\noindent\textbf{Remark.} By the absolute continuity of
the integral, this lemma implies that every curve with finite
$\E_2$ energy satisfies the assumptions of Theorem~\ref{beta-mfd}.

\section{The image of $\Gamma$ is a manifold}

\label{image}

\setnumbers

This section is devoted to the proof of Theorem~\ref{beta-mfd}. We
will argue by contradiction. The proof has two steps; one of them
has preparatory topological character and the second one shows how
to use the assumption on the uniform decay of $\beta$'s.

\medskip
\noindent {\bf Proof of Theorem~\ref{beta-mfd}.}\,\, We recall the
assumption of the theorem that the arclength parametrization
$\Gamma:S_L\to\R^n$ with image $\gamma=\Gamma(S_L)$ satisfies
\eqref{sup-b} for some continuous nondecreasing function
$\omega:[0,L]\to\R$ with $\omega(0)=0.$ In addition, however, we
assume that $\gamma$ is neither homeomorphic to the unit circle
$\S^1$ nor to the unit interval $I=[0,1]$. Our goal is to show
that this leads to a contradiction.

\medskip

\noindent
{\bf Step 1. Triple junctions.}\,\,

\noindent
{\bf Claim:} {\it There exists a triple junction $x\in\gamma$, i.e.
there are three closed sets $\alpha_i\subset\gamma$, $i=1,2,3,$
such that $\alpha_i$ is a continuous image of the unit interval
with $\diam \alpha_i >0$ for $i=1,2,3,$ and such that
\begin{equation} \alpha_i \cap
\alpha_j = \{x\} \qquad\mbox{whenever $i\not=j$, $i,j=1,2,3.$}
\label{junction}
\end{equation}
}

\medskip\noindent\textbf{Remark.} We allow the $\alpha_i$ to have
self-intersections, i.e. we do not require $\alpha_i$ to be a
homeo\-morphic image of the interval. Moreover, more than three
arcs of the curve may meet at $x$; we just need three of them to
obtain the desired contradiction in Step 2  in order to
complete the proof of Theorem~\ref{beta-mfd}.

\medskip\noindent\textbf{Proof of the claim.} We consider two distinct cases.

\smallskip\noindent\textbf{Case 1.} Suppose that $\gamma$ contains
a proper closed subset $\gamma_1$ that is homeomorphic to $\S^1$.
Take a point $y\in \gamma\setminus \gamma_1$, $y=\Gamma(s)$.
Suppose w.l.o.g. that $\Gamma(s_1)\in \gamma_1$ for some $s_1>s$,
$s_1\in [0,L]$ (otherwise just reverse the parametrization). Let
\[
\sigma_0:=\inf\{\sigma>s\colon \Gamma(\sigma)\in \gamma_1\}
\]
It is easy to see that $x=\Gamma(\sigma_0)\in \gamma_1$ is a
triple junction; two of the arcs $\alpha_i$ of $\gamma$ are
contained in $\gamma_1$ and the third one joins $y\not\in\gamma_1$
to $x$.

\smallskip\noindent\textbf{Case 2.} Suppose that Case~1 fails and $\gamma$ contains
no proper closed subset homeo\-morphic to $\S^1$. Consider the
family of all proper subarcs of $\gamma$,
\[
\mathscr{A}=\{\tilde{\gamma}\subset \gamma \colon \
\tilde\gamma\mbox{
is homeomorphic to }I\},
\]
which is partially ordered by inclusion. We will prove in detail
below that every chain in $\A$ has an upper bound in $\A$, so
that  by the Kuratowski--Zorn Lemma $\A$ has a maximal element,
$\gamma_{\textnormal{max}}$. We have
$\gamma_{\textnormal{max}}\not=\gamma$, as $\gamma$ is not
homeo\-morphic to $I$ by assumption. Now, take a point $y\in
\gamma\setminus \gamma_{\textnormal{max}}$, $y=\Gamma(s)$, and
proceed like in Case 1  joining $y$ with an arc to a point
$x\in\gamma_{\textnormal{max}}.$ Notice that $x$ cannot be an
endpoint of $\gamma_{\textnormal{max}}$, since this would
contradict the maximality of $\gamma_{\textnormal{max}}$.

It remains to be shown that every chain in $\A$ indeed has an
upper bound in $\A$, which is obvious for any finite chain. For an
infinite chain $\mathscr{C}:=\{\gamma_l\}_{l\in\Sigma}\subset\A$
where the index may be chosen to coincide with the length of the
respective arc, $l=\mathscr{L}(\gamma_l)\le \mathscr{H}^1(\gamma)$
for $\gamma_l\in\mathscr{C}$, i.e. where the index set $\Sigma$ is
a (in general uncountable) subset of $[0,\mathscr{H}^1(\gamma)]$,
we can choose a nondecreasing sequence of indices $l_i$ with
$\gamma_i\equiv\gamma_{l_i}\in \mathscr{C}$,
$\gamma_i\subset\gamma_{i+1}$ and $l_i\to l^*:=\sup\Sigma\in
(0,\mathscr{H}^1(\gamma)]$.\footnote{Assuming that at least one
member of $\mathscr{C}$ has positive diameter, otherwise the claim
is trivially true.} Now continuously extend the corresponding
nested injective arclength parametrizations
\begin{equation}\label{extension}
\Gamma_i:[-l_i/2,l_i/2]\to\R^n\quad\textnormal{with $\Gamma_i([-l_i/2,l_i/2])
=\gamma_i$\, and
$\Gamma_{i+1}|_{[-l_i/2,l_i/2]}=\Gamma_i$ for all $ i\in\N$}
\end{equation}
by virtue of
$$
\Gamma_i(t):=\begin{cases}
\Gamma_i(-l_i/2)& \Fo t\in [-l^*/2,-l_i/2)\\
\Gamma_i(l_i/2)& \Fo t\in (l_i/2,l^*/2]
\end{cases}
$$
to all of $[-l^*/2,l^*/2].$ Since $|\Gamma_i'(t)|\le 1$ for all
$t\in [l^*/2,l^*/2]$, $i\in\N$, we obtain the uniform bound
$\|\Gamma_i\|_{C^{0,1}([-l^*/2,l^*/2],\R^n)}\le C$ for all $i\in\N,$
which implies by the Theorem of Arzela-Ascoli that some subsequence
$\Gamma_j$ converges to some curve $\Gamma\in C^{0,1}([-l^*/2,l^*/2],\R^n)$
uniformly on $[-l^*/2,l^*/2].$ For distinct parameters $s,t\in (-l^*/2,
l^*/2)$ one can find $j_0\in\N$ such that for all $j\ge j_0$ we have
$s,t\in (-l_j/2,l_j/2),$ so that by \eqref{extension}
$$
|\Gamma(s)-\Gamma(t)|=\lim_{j\to\infty}|\Gamma_j(s)-\Gamma(t)|
\overset{\eqref{extension}}{=}|\Gamma_{j_0}(s)-\Gamma_{j_0}(t )|
\not= 0,
$$
which means that $\Gamma$ is injective, hence a homeomorphism on
the open interval $(-l^*/2,l^*/2)$. But if $\Gamma(l^*/2)$ were
equal to $\Gamma(\tau)$ for some $\tau\in [-l^*/2,l^*/2)$ then the
arc $\Gamma([\tau,l^*/2])$ would be homeomorphic to $\S^1$ which
would contradict our assumption that $\gamma$ is neither
homeomorphic to $\S^1$ nor contains a proper closed subset
homeomorphic to $\S^1$. The same contradiction would occur if
$\Gamma(-l^*/2)=\Gamma(\tau)$ for some $\tau\in (-l^*/2,l^*/2].$
Hence $\gamma^*:=\Gamma([-l^*/2,l^*/2])$ is homeomorphic to the
unit interval $I$, that is $\gamma^*\in \A$. Finally $\gamma^*$ is
maximal for the chain $\mathscr{C}$. Indeed, if
$l^*=\sup\Sigma\in\Sigma$ then $\gamma^*$ is the desired upper
bound because for $l<l^*$ it cannot be that $\gamma_{l^*}$ is
contained in $\gamma_l$, so that total ordering in the chain
implies that $\gamma_l\subset\gamma_{l^*}$. If $l^*\not\in\Sigma$,
on the other hand, we have $l<l^*$ for any $l\in\Sigma,$ which
implies that the corresponding arc $\gamma_l$ is contained in one
of the $\gamma_i$ for $i$ sufficiently large, and hence also
$\gamma_l\subset\gamma^*.$
\heikodetail{

\bigskip

ORIGINAL VERSION WAS: take an arbitrary
arc $\gamma_l\in\mathscr{C}$. If $\gamma_l\subset\gamma_j$ for some
$j\in\N$ then $\gamma_l\subset\gamma^*$ by construction.
So we only have to consider the case that all $\gamma_j\subset\gamma_l$
since $\mathscr{C}$ as a chain is totally ordered. If there was one
point $x\in\gamma_l\setminus \gamma^*$ then we could
estimate
$$
0<\dist(x,\gamma^*)\le\dist(x,\gamma_j)\quad\Foa
j\in\N,
$$
since all arcs under consideration are closed sets and $\gamma_j\subset
\gamma^*
$ for all $j\in\N.$ Hence
$$
\sup\Sigma=l^*\ge\mathscr{L}(\gamma_l)
\ge\dist(x,\gamma_j)+\mathscr{L}(\gamma_j)
\ge\dist(x,\gamma^*)+\mathscr{L}(\gamma_j)
\to \dist(x,\gamma^*)+l^*
$$
as $j\to\infty$,  a contradiction.

\bigskip

{\tt

\bigskip

\raggedright \yy\yy The whole passage above is OK. I corrected
only three really minor things.

\smallskip

There is one  remark to consider: if $l^\ast=\sup\Sigma\in
\Sigma$, then it is practically obvious that $\gamma_{l^\ast}$ is
the desired upper bound. If $l^\ast \not\in \Sigma$, then an
arbitrary arc in our chain is in one of the $\gamma_i$, hence the
last few lines (the second case in your proof of maximality) could
be in fact omitted, modulo a short explanation earlier, to give us
a slightly shorter version.

\smallskip

Am I correct? What do you think?\yy\yy

I AGREE \xx\xx and I have inserted a slightly shorter explanation

\bigskip

}

}

The proof of our claim on the existence of (at least one) triple
junction is complete now.

\noindent
{\bf Step 2. Tilting tubes.}\,\,
We now fix a point $x\in \gamma$ that is a triple junction, and a
small distance $d_0$,
\[
0< d_0< \frac 12 \min_{i=1,2,3} \bigl( \diam \alpha_i\bigr)
\, ,
\]
where $\alpha_i$ denote the closed, connected subsets of $\gamma$
satisfying \eqref{junction} above.

Let $h(s):=s\omega(s)$ for $s\in [0,L]$. Shrinking $d_0$ if
necessary, we can ensure  the initial smallness condition
\begin{equation}\label{init-small}
 h(d_0)< \frac 1{20} d_0.
\end{equation}

Rotating and translating the coordinate system in $\bbbr^n$, we
can assume without loss of generality that $x=0\in \bbbr^n$ and
select the three distinct points
\[
y_i \in \alpha_i \cap \partial B(0,d_0), \qquad i=1,2,3
\]
where $y_1=(d_0,0,\ldots, 0)$. Assumption \eqref{sup-b} implies
now
\begin{equation}
\label{box-0} \gamma \cap B(0,d_0) \ \subset \ U_{2 h(d_0)}
(G(x,y_1))\, .
\end{equation}
\heikodetail{

\bigskip

There exists an (infimal) line $G_x$ through $x=0$ such that
$$
\sup_{y\in\gamma\cap B(0,d_0)}\dist(y,G_x)\le d_0\omega(d_0)=h(d_0).
$$
Now (see my old drawing ad 3.2 -1-)
\begin{eqnarray*}
\dist(y,G(x,y_1)) & \le & |y-\pi_{G_x}(y)|+\dist(\pi_{G_x}(y),G(x,y_1))\\
& = & |y-\pi_{G_x}(y)|+|\pi_{G_x}(y)|\cdot\sin\ang (G_x,G(x,y_1))\\
& \le & h(d_0)+|y|\cdot\sin\ang (G_x,G(x,y_1))\\
& \le & h(d_0)+d_0\cdot\sin\ang (G_x,G(x,y_1))\\
& \le & h(d_0)+d_0\omega(d_0)=2h(d_0),
\end{eqnarray*}
since $\dist(y_1,G_x)\le d_0\omega(d_0)$ as $y_1\in B(0,d_0)$, and
$$
\dist(y_1,G_x)=|y_1-\pi_{G_x}(y_1)|=|y_1|\cdot\sin\ang(G_x,G(x,y_1)).
$$

\bigskip

}%
The intersection of the sphere $\partial B(0,d_0)$ with the tube
$ U_{2 h(d_0)} (G (x,y_1))$ consists of two symmetric spherical
caps; by Dirichlet's pigeon-hole principle, one of these caps must
contain two of the three distinct points $y_i$. Renumbering the
$\alpha_i$ and $y_i$ if necessary, we may assume that $y_1$ is as
above and $y_2=(a,y_2')\in\alpha_2\cap\partial B(0,d_0) $ with
$a>0$ and $y_2'\in \bbbr^{n-1}$, $|y_2'|\le 2 h(d_0)$.

Let $v_0=(-1,0,\ldots, 0)$ and $H_0=(v_0)^\perp$. Fix a point
$z\in \alpha_1\cap \bigl(\frac 12 y_1 + H_0\bigr)$.

From now on, we will work only with $\alpha_1$ and
$\alpha_2$. Proceeding inductively, we shall define a
sequence of distances $d_m\to 0$, unit vectors $v_m$, linear
$(n-1)$-dimensional subspaces $H_m=(v_m)^\perp$ and points $x_m\in
\alpha_2$ such that
\begin{equation}
\label{dist-zxm} |z-x_m|\le 2 h(d_m), \qquad m=1,2,\ldots
\end{equation}
As $d_m\to 0$ and $h(s)\to 0$ as $s\to 0$, this will yield $z=\lim
x_m \in \alpha_1 \cap \alpha_2$, a contradiction.

The distances $d_m$, auxiliary vectors $v_m$ and hyperplanes
$H_m=(v_m)^\perp$ will be defined in such a way that for all
$m=1,2,\ldots$
\begin{gather}
4h_{m-1} \le d_m \le 6h_{m-1} \qquad \mbox{where $h_m:= h(d_m)$,}\label{dmm-1} \\
\ang(v_m,v_{m-1}) \le \frac \pi 4, \label{vmm-1}\\
z_m= z+ d_m v_m \in \alpha_1, \label{choice-zm} \\
\gamma \cap B(z,d_m) \subset  U_{2 h_m} (G_m) \qquad \mbox{where
$G_m=G(z,z_m)$.}\label{box-m}
\end{gather}
For $P_m(t)= z+tv_m + H_m$ we shall also show that
\begin{equation}
\label{cut-m} P_m(t) \cap \alpha_i \cap U_{2
h_m}(G_m)\not=\emptyset \qquad\mbox{for all $|t| \le \frac 12 d_m$
and $i=1,2$,}
\end{equation}
for each $m=1,2,\ldots$. Notice that \eqref{dmm-1} in connection
with the initial smallness condition \eqref{init-small} will yield
$d_m\to 0$ as $m\to\infty.$

We begin the construction for $m=1$. Select $z_1\in P_0(4 h_0)
\cap \alpha_1$, $h_0=h(d_0)$. Such a point exists since $\alpha_1$
joins $z$ to $x=0$ and by continuity must intersect all planes
$z+tv_0 + H_0$, $|t|\le \frac 12 d_0$, while staying in the tube
$U_{2 h(d_0)}(G(x,y_1))$. Let $v_1=(z_1-z)/|z_1-z|$,
$H_1:=(v_1)^\perp$, and $P_1(t):=z+tv_1+H_1$. Note that
$\ang(v_1,v_0)\le \pi/ 4$ by construction. Set $d_1=|z_1-z|$.

We already have \eqref{dmm-1}--\eqref{choice-zm} for $m=1$;
condition \eqref{box-m} for $m=1$ follows directly from
\eqref{sup-b}. To obtain \eqref{cut-m} for $m=1$, we just use
\eqref{vmm-1} and continuity.

Assume now that $d_m$, $v_m$, $H_m$,  $z_m$, and $P_m$ have
already been defined for $m=1,\ldots, N$ so that
\eqref{dmm-1}--\eqref{cut-m} are satisfied for all $1\le m\le N$.
We use \eqref{cut-m} for $m=N$ to select a point $z_{N+1}$,
\[
z_{N+1} \in U_{2 h_N}(G_N) \cap P_N(2h_N) \cap \alpha_1\, .
\]
Clearly, $4 h_N \le |z_{N+1}-z|\le 6 h_N$ (the second estimate is
a simple application of the triangle inequality). Thus,
$d_{N+1}:=|z_{N+1}-z|$ satisfies \eqref{dmm-1} for $m=N+1$, and
choosing $v_{N+1}:=(z_{N+1}-z)/|z_{N+1}-z|$ we also have
\eqref{vmm-1}--\eqref{choice-zm} for $m=N+1$.

Again, \eqref{box-m} for $m=N+1$ follows from the assumption on
the decay of $\beta$'s. Thus, the intersection $\alpha_i \cap
B(z,d_{N+1}) \subset U_{2 h_{N+1}}(G_{N+1})$, $i=1,2$; combining
these inclusions with \eqref{vmm-1} and with continuity, we obtain
\eqref{cut-m} for $m=N+1$.

This completes the inductive construction. Now, using
\eqref{cut-m} for $i=2$, we select for each $m$ a point
\[
x_m \in U_{2h_m}(G_m) \cap P_m(0) \cap \alpha_2\, .
\]
By definition of $ U_{2h_m}(G_m)$, \eqref{dist-zxm} does hold.
This completes the whole proof of Theorem~\ref{beta-mfd}.

\heikodetail{

\bigskip

For all these calculations see my notes on ad 3.2 -2- and -3-
and the drawing there. It might be worth thinking of a
drawing as on ad3.2-2-...?\xx One important aspect of the
geometry is that e.g. the (even) plane $P_1(d_0) $ does not
intersect the left halfspace in the slab , i.e.
$$
\{\xi\in\R^n: \xi_1<0\}\cap U_{2h_0}(G(x,y_1))\cap P_1(d_0)=\emptyset
$$
which implies that $\alpha_3\cap B(0,d_0)$ does not come into
play, so that the continuity argument really yields points in
$\alpha_1$ and $\alpha_2.$ Should we include this somewhere
above??\xx\xx{} \yy I added one phrase.

\bigskip

}

\section{Differentiability}

\label{diff}

\setnumbers

\heikodetail{

\yy\yy I have adjusted this Section slightly, for the sake of the
application to knots. One remark: I think that the old extra
restriction that the distance of the two points on $\gamma$ be
${}< \frac 12 \diam \gamma$ is not really necessary (in other
words: there is an implicit relation between the diameter and
energy of $\gamma$). In fact we do not seem to use this
restriction anywhere. If the points are too far away from each
other, then it may  simply turn out that the balance condition
\eqref{ed} is be violated.

\smallskip

Am I right?  \yy\yy

\xx\xx In the few lines of the proof of Prop. \ref{prop:5.1} we
need that $\partial B (x,d)\cap\gamma$ is nonempty, that is probably
the reason we introduced $\diam\gamma/2$ as an additional bound for
$d$. I am not sure if we can get rid of it...

\yy\yy Yes, we need to know that $\partial B (x,r)\cap\gamma$ is
nonempty for $r=d/2^N$. But once we have a couple $x,y$ such that
$|x-y|=d$, then we know that the curve must go from $x$ to $y$,
and all the spheres $\partial B (x,r)$, $r\in (0,d)$, are
intersected by $\gamma$. This justifies Step 1 in the proof below,
doesn't it?

If I remember correctly, the bound ${}< \frac 12 \diam \gamma$ was
necessary since at some place we wanted to know that the curve was
not totally contained in $B(x,d)$. As far as I can see we do not
use this here (or in Section 5) --- injectivity is now taken for
granted.

\xx\xx I think you are correct, and we are on the safe side here.\xx\xx

\bigskip

}

Throughout this section, we fix $q>2$ and consider a rectifiable
curve $\gamma=\Gamma(S_L)$ whose arc\-length parametrization
$\Gamma$ is injective on $S_L$. The first step towards the proof
of Theorem~\ref{thm:1.3} is to establish the following.

\begin{proposition}
\label{prop:5.1} Let $q>2$. Assume that $\Gamma\colon S_L\to
\bbbr^n$ is injective and $\E_q(\Gamma)< E < \infty$. Then
$\Gamma'$ is well defined everywhere and $\Gamma'\in C^\kappa$ for
$\kappa:=\frac{q-2}{q+4}\in (0,1)$.

Moreover there exist two positive constants $\delta_2(q)$,
$c_2(q)$ such that whenever $x=\Gamma(s)$ and $y=\Gamma(t)$
satisfy $|x-y|=d<\delta_2(q) E^{-1/(q-2)} $, then
\begin{equation}
\label{phi-angle} \phi: = c_2(q) E^{1/(q+4)}d^{(q-2)/(q+4)} <
\frac 14
\end{equation}
and we have
\begin{gather}
\label{osc-tan} |\Gamma'(s)-\Gamma'(t)|\le {c_2(q)}\,
E^{1/(q+4)}
|\Gamma(s)-\Gamma(t)|^{\kappa}, \\
\label{bilip}
\frac 34 |s-t| \le |\Gamma(s)-\Gamma(t)|\le |s-t|, \\
\label{2cones} \gamma \cap B(x,2d) \cap B(y,2 d) \ \subset \
C_\phi(x,y) \cap C_\phi(y,x).
\end{gather}
\end{proposition}

\medskip\noindent\textbf{Proof.} The argument is in fact similar to the
proof of Corollary~2.6 and Theorem~2.10 in \cite{ssvdm-triple}. We
just sketch the main points, leaving (relatively easy)
computational details as an exercise.

Fix $x,y\in \gamma$ with $0< |x-y|= d.$

\smallskip\noindent\textbf{Step 1.}  For $N=0,1,2\ldots$ set $d_N=d/2^N$, and
select points $y_N\in \partial B(x,d_N)\cap \gamma$ so that
$y_0=y$. Let
\begin{equation}
\eps_N:=  \bigl(c_0(q) E\bigr)^{1/(q+4)} d_N{}^{\kappa}
\label{eps-N}
\end{equation}
so that condition \eqref{ed} of Lemma~\ref{beta} is satisfied for
$\eps_N$ and $d_N$. The lemma yields
\begin{equation}
\gamma\cap B(x,2d_N) \subset U_{20\eps_Nd_N}(G(x,y_N)), \qquad
N=0,1,2,\ldots
\end{equation}
so that the lines $G_N:=G(x,y_N)$ satisfy
\begin{equation}
\sin\ang (G_N, G_{N+1}) \le \frac{20\eps_Nd_N}{d_{N+1}} =
40\eps_N\, .
\end{equation}
Thus, $\phi_N:=\ang(G_N,G_{N+1})\le 80 \eps_N$. Using
\eqref{eps-N} and summing a geometric series (here the assumption
$q>2$ is crucial!), we obtain
\begin{equation}
\label{def-phi}  \sum_{N=0}^\infty \phi_N \le \phi:= c_2(q)
E^{1/(q+4)}d^{\kappa}
\end{equation}
where $c_2(q)=80\,  c_0(q)^{1/(q+4)}\sum_{N=0}^\infty
2^{-N\kappa}$. Now, to guarantee $\phi<1/4$, one just assumes that
$d$ is sufficiently small, i.e. $d< \delta_2(q)E^{-1/(q-2)}$ with
$\delta_2(q):=(4c_2(q))^{-1/\kappa}$.  By induction,
\begin{equation}
\gamma\cap B(x,2d) \ \subset\  C_{\phi_0+\cdots+\phi_N} (x,y)
\cup (U_{20\eps_Nd_N}(G_N)\cap B(x,2d_N))\, .
\end{equation}
Passing to the limit $N\to \infty$, we obtain
\begin{equation}
\gamma \cap B(x,2d) \subset C_\phi (x,y)
\end{equation}
with $\phi\equiv \phi(q,E,d)$ defined by \eqref{def-phi}.

\smallskip\noindent\textbf{Step 2.} Reversing the roles of $x$
and $y$ we obtain
\[
\gamma \cap B(x,2d) \cap B(y,2 d) \ \subset \ C_\phi(x,y) \cap
C_\phi(y,x)
\]
where $\phi$ is defined by \eqref{def-phi}; this is the desired
condition \eqref{2cones}.

\smallskip\noindent\textbf{Step 3.} Assume now that
%$x=\Gamma(s)$,
%$y=\Gamma(t)$, where
%\[
%|s-t|\le \delta_2(q) E^{-1/(q-2)}\, .
%\]
%Then $|x-y|=d$ satisfies the constraints needed in Steps 1 and~2.
%If
$\Gamma$ is differentiable at $s$ and $t$ and recall that $\Gamma$
was supposed to be injective. Condition \eqref{2cones} yields then
\begin{equation}
\label{ang-st} \ang(\Gamma'(s),\Gamma'(t))\le \phi = c_2(q)
E^{1/(q+4)}d^{\kappa} = c_2(q)
E^{1/(q+4)}|\Gamma(s)-\Gamma(t)|^{\kappa}\,
\end{equation}
(note that the difference quotients of $\Gamma$ at $s$ and $t$
must belong to cones with vertices at $0$, axis parallel to $y-x$
and opening angle given by \eqref{phi-angle}).

\smallskip\noindent\textbf{Step 4.} Since $\Gamma$ is
differentiable everywhere, and $|\Gamma'|=1$ a.e., \eqref{ang-st}
gives \eqref{osc-tan} on a (dense) set of full measure. Thus,
$\Gamma'$ has a continuous extension $F$ to all of $S_L$;
one easily checks that in fact $F =\Gamma'$
\emph{everywhere\/}. Finally, assuming without loss of
generality that $t>s$, we estimate
\begin{eqnarray*}
|\Gamma(t)-\Gamma(s)| & \ge & \langle
\Gamma(t)-\Gamma(s),\Gamma'(s)\rangle\\ & = & \left\langle
\int_s^t \bigl(\Gamma'(\tau)-\Gamma'(s)+\Gamma'(s)\bigr)\,
d\tau,\,
\Gamma'(s)\right\rangle \\
& \ge & (t-s) \left(1 - \sup_{\tau\in [s,t]}
|\Gamma'(\tau)-\Gamma'(s)|\right)  \ \ge \ \frac 34 (t-s)\, .
\end{eqnarray*}
(To check the last inequality, let $S$ be the closed slab bounded
by two planes passing through $x$ and $y$, and perpendicular to
$x-y$, i.e. to the common axis of the two cones, and note that for
each $\tau \in [s,t]$ we have in fact $\Gamma(\tau)\in C_\phi(x,y)
\cap C_\phi(y,x)\cap S$. This follows from the bound
$|\Gamma'(s)-\Gamma'(t)|<1/4$, injectivity of $\Gamma$ and
\eqref{2cones}. Thus, for each such $\tau$ we also have
$|\Gamma'(\tau)-\Gamma'(s)|<1/4$.) The bi-Lipschitz condition
\eqref{bilip} follows.

\smallskip

The proof of Proposition~\ref{prop:5.1} is complete now. (See also
\cite[Proof of Thm. 2.10]{ssvdm-triple} where a similar scheme of
reasoning is used.) \hfill $\Box$

\section{Energy bounds and knot classes}\label{sec:5}

We start this section with the observation that $\E_q$
is repulsive (or charge), that is, $\E_q$ blows up on a sequence
of knots converging uniformly to a limit curve with
self-crossings.
\begin{proposition}\label{self-repulsive}
If $\Gamma:S_L\to\R^n$ is a closed arclength parametrized curve of length
$0<L<\infty$ with $\Gamma(s)=\Gamma(t)$ for different arclength
parameters $s\not= t,$ $s,t\in S_L,$ and if there is a sequence
of rectifiable closed injective curves $\gamma_k:S_L\to\R^n$ converging
uniformly to $\Gamma$, then $\E_q(\gamma_k)\to\infty$ as
$k\to\infty$ for any $q>2.$
\end{proposition}
\proof
Assume to the contrary that (for a suitable subsequence)
$\lim_{k\to\infty}\E_q(\gamma_k)< E<\infty.$ We set
\begin{equation}\label{eps}
\eps:=\frac 12 \min\Big\{\diam\Gamma([s,t]),\diam\Gamma(S_L\setminus [s,t]),
\delta_2(q)E^{\frac{-1}{q-2}}\Big\}>0,
\end{equation}
where $\delta_2(q)$ is the constant of Proposition \ref{prop:5.1},
and choose $\tau\in (s,t)$ and $\sigma\in S_L\setminus [s,t]$ such that
$$
|\Gamma(\tau)-\Gamma(t)| =  \frac 12 \diam\Gamma([s,t])\quad\AND\quad
|\Gamma(\sigma)-\Gamma(s)|=  \frac 12 \diam\Gamma(S_L\setminus [s,t]).
$$
For sufficiently large $k_0=k_0(\eps)\in\N$ we find
$\|\gamma_k-\Gamma\|_{C^0(S_L,\R^n)}<\eps/10$ for all $k\ge k_0.$
In particular, by \eqref{eps},
\begin{eqnarray}
|\gamma_k(\tau)-\gamma_k(t)|&\ge &|\Gamma(\tau)-\Gamma(t)|-\frac{2\eps}{10}=
\frac 12 \diam\Gamma([s,t])-\frac \eps 5\overset{\eqref{eps}}{\ge}\frac 45
\eps,
\notag\\
\label{lb}
\\
\textnormal{and, analagously,\,\,} &&|\gamma_k(\sigma)-\gamma_k(s)|\ge
\frac 45 \eps,
\notag
\end{eqnarray}
but
$$
\delta_k:=|\gamma_k(t)-\gamma_k(s)|\le\frac \eps 5\overset{\eqref{eps}}{<}
\delta_2(q)E^{\frac{-1}{q-2}}\quad\Foa k\ge k_0.
$$
Hence, we can apply \eqref{2cones} of Proposition \ref{prop:5.1} to
obtain
the inclusion
$$
\gamma_k\cap B(\gamma_k(t),2\delta_k)\cap
B(\gamma_k(s),2\delta_k)\overset{\eqref{2cones}}{\subset}
C_\phi(\gamma_k(t),\gamma_k(s))\cap C_\phi(\gamma_k(s),\gamma_k(t)).
$$
Since there is an integer $k_1\ge k_0$ such that
$\E_q<E$ for all $k\ge k_1$ we know that the corresponding
injective arclength parametrizations $\Gamma_k$ are continuously
differentiable according to Proposition \ref{prop:5.1}, so that
the points $\gamma_k(t)$ and $\gamma_k(s)$ must be connected by
a subarc of $\gamma_k$ that is completely contained in the doubly
conical region
$$
D_k:=C_\phi(\gamma_k(t),\gamma_k(s))\cap
C_\phi(\gamma_k(s),\gamma_k(t))\cap B(\frac 12 (\gamma_k(t)+\gamma_k(s)),
\frac{\delta_k}{2})
$$
of diameter $\delta_k\le\eps/5.$ (Otherwise, the unit
tangent vector of the arclength parametrization $\Gamma_k$ would
jump at $\gamma_k(t)$ and $\gamma_k(s)$
contradicting $C^1$-smoothness for  $k\ge k_1$.)
Since all $\gamma_k$ are simple, either the point $\gamma_k(\tau)$, or
$\gamma_k(\sigma)$ lies on that connecting arc within $D_k$, thus
contradicting the lower bound $4\eps/5$ in \eqref{lb}.
\qed

\begin{proposition}
\label{prop:strong} If $q>2$, then the $\E_q$-energy is
\emph{strong\/} in the following sense:  For each $E>0$ and $L>0$
there are at most finitely many knot types which have a
representative $\gamma$ such that
\[
\E_q(\gamma)<E, \qquad \H^1(\gamma)=L.
\]
\end{proposition}

\noindent\textbf{Remark.} The length constraint $\H^1(\gamma)=L$
is necessary here, since by rescaling an arbitrary smooth simple
curve we can make its $\E_q$-energy as small as one wishes. An
alternative would be to consider
$\tilde{\E}_q(\gamma):=(\H^1(\gamma))^{q-2}\E_q(\Gamma)$. This is
a scale invariant energy.

\medskip\noindent
\textbf{Proof.} We argue by contradiction. Assume there are
infinitely many knot types of length $L$ with the same energy
bound, and by translational invariance we can assume moreover that
all these knots contain the origin. Take their arclength
representatives $\Gamma_j$, $j=1,2,\ldots,$ and use inequality
\eqref{osc-tan} of Proposition ~\ref{prop:5.1} to conclude that
the family
\[
\{\Gamma_j'\}_{j=1,2,\ldots} \, \subset \ C^0(S_L,\S^{2})
\]
is eqicontinuous. Invoking the Arzela--Ascoli compactness theorem
and passing to a subsequence, we may assume that $\Gamma_j$
converges in the $C^1$-topology to some limit $\Gamma\in
C^1(S_L,\R^3)$. Let $\gamma$ be the curve parametrized by
$\Gamma$.

We next check that $\gamma$ is simple, i.e. $\Gamma$ is injective
on $S_L\equiv \R/L\Z$. To this end, we shall rely on
Proposition~\ref{prop:5.1} to prove that there exists an
$\eps_0=\eps_0(q,E)>0$ such that all $\Gamma_j$ satisfy
\begin{equation}
\label{low-bound} |\Gamma_j(s)-\Gamma_j(t)|\ge \min \Bigl( \eps_0,
\frac{|s-t|}{2} \Bigr) \qquad\mbox{for all $j$ and all $s,t\in
S_L$.}
\end{equation}
Upon passing to the limit $j\to\infty$, this implies the
injectivity of $\Gamma$. All $\gamma_j$ with $j$ sufficiently
large are contained in a small $C^1$ neighbour\-hood of $\gamma$.
Thus, according to a known isotopy result, see e.g. \cite[Chapter
8]{hirsch} or \cite{blatt}, they would all have to be of the same
knot type, thereby contradicting the assumption that each
$\gamma_j$ is in a different isotopy class.

To complete the proof, it is now enough to prove
\eqref{low-bound}. Consider $g_j\in C^1 (S_L\times S_L) $ given by
\[
g_j(s,t) \colon = |\Gamma_j(s)-\Gamma_j(t)|^2\, .
\]
By Proposition~\ref{prop:5.1} the $\Gamma_j$ are uniformly bounded
in $C^{1,\kappa}$, where $\kappa=(q-2)/(q+4)$. Thus, it is easy to
show that there is a constant $\eps_1=\eps_1(q,E)>0$ such that
\begin{equation}
\label{gj-below} g_j(s,t) \ge \frac{|s-t|^2}4 \qquad\mbox{for all
$j$ and all $s,t$ such that $|s-t| \le \eps_1(q,E)$.}
\end{equation}
Since $\Sigma=S_L\times S_L\setminus \{(s,t)\colon
|s-t|<\eps_1(q,E)\}$ is compact, for each $j\, $ there is  a pair
$(s_j,t_j)\in \Sigma$ such that
\[
g_j(s_j,t_j)\le g_j(s,t) \qquad \mbox{for all $(s,t)\in \Sigma$.}
\]
Now, we either have $|s_j-t_j| = \eps_1(q,E)$ in which case
\eqref{gj-below} implies
\begin{equation}
\label{heiko-4} g_j(s,t) \ge \frac{\eps_1(q,E)^2}4 \qquad\mbox{for
all $s,t\in \Sigma$,}
\end{equation}
or, by minimality, we have $\nabla g_j(s_j,t_j)=0$, which is
equivalent to
\begin{equation}
\label{perp} \Gamma_j'(s_j) \perp
\bigl(\Gamma_j(s_j)-\Gamma_j(t_j)\bigr) \qquad \mbox{and}\qquad
\Gamma_j'(t_j) \perp \bigl(\Gamma_j(s_j)-\Gamma_j(t_j)\bigr).
\end{equation}
Fix $j$. Let $d_j:=|\Gamma_j(s_j)-\Gamma_j(t_j)|$. If $ d_j <
\delta_2(q) E^{-1/{q-2}}, $ where $\delta_2(q)$ stands for the
constant from Proposition~\ref{prop:5.1}, then, by
\eqref{phi-angle} and \eqref{2cones} of that Proposition, we have
\begin{equation*}
\label{phi-angle-j} \phi_j: = c_2(q) E^{1/(q+4)}d_j^{\kappa} <
\frac 14
\end{equation*}
and
\[
\gamma_j \cap B(\Gamma_j(s_j),2d_j) \cap B(\Gamma_j(t_j),2d_j) \
\subset \ C_{1/4}(\Gamma_j(s_j),\Gamma_j(t_j)) \cap
C_{1/4}(\Gamma_j(t_j),\Gamma_j(s_j)).
\]
The last condition, however, clearly contradicts \eqref{perp}.
Hence,
\begin{eqnarray}
d_j= |\Gamma_j(s_j)-\Gamma_j(t_j)| & = & \inf_{(s,t)\in \Sigma}
|\Gamma_j(s)-\Gamma_j(t)| \nonumber \\
& \ge &  \eps_2(q,E)\, \colon=\, \delta_2(q) E^{-1/{q-2}}
\qquad\mbox{for each $j=1,2,\ldots$} \label{dj-below}
\end{eqnarray}
Summarizing \eqref{gj-below}, \eqref{heiko-4} and
\eqref{dj-below}, we obtain \eqref{low-bound} with $\eps_0\colon =
\min \bigl\{\eps_1(q,E)/2, \eps_2(q,E)\bigr\}$.  \hfill $\Box$

\medskip

%\noindent \textbf{Left:} The plane which passes through two points
%$\Gamma_j(s_j),\Gamma_j(t_j)$, and contains $\Gamma'_j(s_j)$.
% The planar cross-section of
%\[
%C_{\pi /4}\bigl(\Gamma_j(s_j); \Gamma_j(t_j)\bigr)\ \cap\ C_{\pi
%/4}\bigl(\Gamma_j(t_j); \Gamma_j(s_j)\bigr)
%\]
%is shaded. The tangent vectors to $\gamma_j$ at these two points
%are perpendicular to the common axis of the cones.

\medskip

Now
we present the proof of the isotopy result,
Theorem~\ref{thm:1.2}. The proof consists of two steps. The first
one, see Proposition~\ref{poly-1} below, is preparatory: we use
Proposition~\ref{prop:5.1} to show that a curve $\gamma$ of length
$L$ and finite energy at most $E$ is ambient isotopic to a
polygonal line that has roughly $LE^{1/(q-2)}$ vertices, all of
them belonging to $\gamma$. In the second step, we replace two
curves that are close in Hausdorff distance by polygonal curves
(staying in the same knot class)  and exhibit a series of $\Delta$
and $\Delta^{-1}$-moves\footnote{These are \emph{not} the
so-called Reidemeister moves; see \cite[Chapter~1]
{burde-zieschang} for the distinction.} transforming one of them
into the other one. (The proof that we present gives a value of
$\delta_3$ which is far from optimal; we do not know how to obtain
a sharp result of that type.)

Before passing to the details, let us recall a definition, see
e.g. \cite[Chapter~1]{burde-zieschang}.

\begin{definition}\rm
Let $u$ be one of the segments of a polygonal knot $\gamma$ in
$\R^3$ and let $T=\mathrm{conv}(u,v,w)$ be a triangular surface
bounded by the segments $u,v,w$ such that $T\cap \gamma= u$. We
say that
\[
\gamma' = (\gamma\setminus u) \cup v \cup w
\]
\emph{results from $\gamma$ by a $\Delta$-move\/}. The inverse
operation is called a \emph{$\Delta^{-1}$-move}.
\end{definition}

Let $\gamma_1$ and $\gamma_2$ be two polygonal knots in $\R^3.$ If
$\gamma_1$ can be obtained from $\gamma_2$ by a finite sequence of
$\Delta$ and $\Delta^{-1}$-moves, then one says that $\gamma_1$
and $\gamma_2$ are \emph{combinatorially equivalent\/}. Two
polygonal knots $\gamma_1$ and $\gamma_2$ are ambient isotopic if
and only if they are combinatorially equivalent, see
\cite[Chapter~1]{burde-zieschang}.

\begin{proposition}\label{poly-1}
Let $q>2$. Assume that $\Gamma\colon S_L\to \R^3$ is injective and
$\E_q(\Gamma)<E$. Let $\delta_2(q)>0$ be the constant defined in
Proposition~\ref{prop:5.1}. Then $\gamma=\Gamma(S_L)$ is ambient
isotopic to the polygonal curve
\[
P_\gamma = \bigcup_{i=1}^N [x_i,x_{i+1}]
\]
with $N$ vertices $x_i=\Gamma(t_i)\in \gamma$, whenever the
parameters $0= t_1<\ldots < t_N< L$ and $t_{N+1}=t_1$ are chosen
on $S_L$ so that
\begin{equation}\label{spacing}
|x_i-x_{i+1}|< \delta_2(q) E^{-1/(q-2)}\, .
\end{equation}
\end{proposition}

\medskip\noindent\textbf{Proof.} We follow \cite[Prop.~5.2]{marta}
with minor technical changes. For $x\not=y\in \R^3$ we denote the
closed halfspace
\[
H^+(x,y) \colon= \{z\in \R^3\colon \langle z-x, y-x\rangle \ge
0\}\, .
\]
We shall work with `double cones'
\[
K(x,y)\colon= C_{1/4}(x,y) \cap C_{1/4}(y,x) \cap H^+(x,y) \cap
H^+(y,x)\, .
\]
\heikodetail{\tt\raggedright

\xx I have tried to find the exact term for that geometric object
but was not successful, but there must be a proper name for that.
Our diamond-phrase is a little misleading since the reader would
expect corners...

\yy I google'd the term `double cone' and, among other things, saw
a film showing a well-known tricky physical experiment where a
wooden double cone seems to `roll uphill' (but, obviously, the
center of gravity moves in the correct direction). So, at least we
are not the only ones...

}
%Fix $0=t_0<t_1<\ldots < t_N< L$, let $t_{N+1}:=t_1$ and
%$x_i:=\Gamma(t_i)$ for $i=1,\ldots,N+1$. Assume that the $t_i$'s
%are chosen so that the spacing condition \eqref{spacing} holds.
For sake of brevity, set $K_i:=K(x_i,x_{i+1})$ and
$v_i:=x_{i+1}-x_{i}$. We are going to use
Proposition~\ref{prop:5.1} to verify two properties of $K_i$.

\smallskip\noindent\textbf{Claim 1.} \emph{For each $z\in K_i$
the intersection of $\gamma$ and the two-dimensional disk
\[
D_i(z):= K_i\cap (z+v_i^\perp)
\]
contains precisely one point. If $\diam D_i(z)>0$, then this point
of $\gamma$ is in the interior of $K_i$.}

\smallskip

Indeed, note first that $\gamma \cap D_i(z)$ is nonempty, as an
arc of $\gamma$ joining $x_i$ with $x_{i+1}$ must be contained in
$K_i$ since if this were not the case, then \eqref{2cones} of
Proposition~\ref{prop:5.1} would be impossible for an injective
and differentiable $\Gamma$. If there were two distinct points
$y_1,y_2\in \gamma \cap D_i(z)$, then \eqref{2cones} could not
hold both for the couple $x=x_i,y=y_1$, and for the couple $x=x_i,
y=y_2$, simultaneously. Finally, the second statement of Claim~1
follows from the fact that Inequality \eqref{phi-angle} is strict.

\smallskip\noindent\textbf{Claim 2.} \emph{Whenever $i\not=j \pmod N$ %\textnormal{mod$(N)$}$,
we find that the sets $K_i\setminus\{x_i,x_{i+1}\}$ and
$K_j\setminus\{x_j,x_{j+1}\}$ are disjoint. }
\heikodetail{\raggedright\tt

\xx Since I did not quite understand the last lines of your
argument I thought of an alternative, what do you think of the
following argument?

\yy It is OK, and easier to follow than the old one. I corrected
just the spot in (5.16) and two indices below. Is it all right
now?
}

Suppose to the contrary that
\begin{equation}\label{A}
(K_i\setminus\{x_i,x_{i+1}\})\cap (K_j\setminus\{x_j,x_{j+1}\})\not=\emptyset,
\end{equation}
and assume without loss of generality
\begin{equation}\label{B}
\diam K_j\le\diam K_i.
\end{equation}
If $x_j=\Gamma(t_j)$ were contained in $K_i\setminus\{x_i,x_{i+1}\}$ then
we would either find that the disk $D_i(x_j)$ contains two distinct curve
points contradicting Claim 1, or that there is a parameter $\tau\in (t_i,t_{i+1})$
such that $\Gamma(\tau)=\Gamma(t_j)$ although $\Gamma$ is injective, which is
absurd.
The same reasoning can be applied to $x_{j+1}=\Gamma(t_{j+1}),$ so that we conclude
from \eqref{2cones} and Assumptions \eqref{A} and \eqref{B} that the two
tips
$x_j,$ $x_{j+1}$ of $K_j$ are contained in the set $Z_i$ defined as
\begin{equation}\label{C}
Z_i:=C_{\frac 14}(x_i,x_{i+1})\cap C_{\frac 14} (x_{i+1},x_i)\cap
B(x_i,2|v_i|)\cap B(x_{i+1},2|v_i|)\setminus \Big[K_i\setminus\{x_i,x_{i+1}\}\Big],
\end{equation}
which is just the intersection of the two cones within the balls centered
in $x_i$ and $x_{i+1}$ but without the slab bounded by the two parallel
planes $\partial H^+(x_i,x_{i+1})$ and
$\partial H^+(x_{i+1},x_i)$.

We know that $x_j\not=x_i$ since $i\not= j\pmod N$ and $\Gamma$ is
injective. If $x_j\not=x_{i+1}$ then \eqref{A}, \eqref{B}, and
\eqref{C}  enforce
$$
|v_i|\overset{\eqref{B}}{\ge}|v_j|\overset{\eqref{A}}{>}\min\{|x_j-x_{i+1}|,|x_j-x_i|\},
$$
and
\begin{equation}\label{DD}
x_{j+1}\in\textnormal{int}(H^+(x_i,x_{i+1}))\,\,\cap\,\,
\textnormal{int}(H^+(x_{i+1},x_i)),
\end{equation}
which by \eqref{2cones} leads to $x_{j+1}\in K_i$ contradicting
\eqref{C} unless $x_{j+1}=x_i$. If in the latter case $x_{j}$
is contained in $\R^3\setminus H^+(x_{i+1},x_i)$ then we obtain
$|v_j|=|x_{j+1}-x_j|>|v_i|$ contradicting our assumption
\eqref{B}. If, on the other hand, $x_{j}$ is in
$H^+(x_{i+1},x_i)$, it is by \eqref{C} actually contained in
$\R^3\setminus H^+(x_i,x_{i+1})$, but then \eqref{A} cannot hold.

Finally, $x_j=x_{i+1}$ in combination with \eqref{A} also leading to
\eqref{DD} is a contradictory statement, since $|v_j|\le |v_i|$ by \eqref{B}.

\heikodetail{

\bigskip

Pavel's argument:

\smallskip

Suppose the contrary. Let $i<j$ and assume $z\in \mathrm{int}\,
K_i \cap \mathrm{int}\, K_j$. Then $\diam D_i(z) \le \tan \frac 18
|x_{i+1}-x_i| < \frac 12 \delta_2(q) E^{-1/(q-2)}$. Similarly,
$\diam D_j(z) < \frac 12 \delta_2(q) E^{-1/(q-2)}$. Thus, using
Claim~1 above, we can find two parameters, $s_i\in (t_i,t_{i+1})$
and $s_j\in (t_j,t_{j+1})$ such that $\Gamma(s_l)\in D_l(z)$,
$l=i,j$. Hence,
\[
|\Gamma(s_i)-\Gamma(s_j)| < |\Gamma(s_i)-z|+|z-\Gamma(s_j)|<
\delta_2(q) E^{-1/(q-2)}\, .
\]
Thus, again by Proposition~\ref{prop:5.1} and the injectivity of
$\Gamma$, the image of one of the arcs of $S_L$ between $s_i$ and
$s_{j+1}$ is connected in $K:=K(\Gamma(s_i),\Gamma(s_j))$. Without
loss of generality, we assume that this image contains
$x_{i+1}=\Gamma(t_{i+1})$ and $x_j=\Gamma(t_j)$. Now, since
$\Gamma(s_i)\in K_i$ and $\Gamma(s_j)\in K_j$, all four angles
\begin{gather*}
\ang\bigl(v_i,x_{i+1}-\Gamma(s_i)\bigr), \qquad
\ang\bigl(x_{i+1}-\Gamma(s_i),\Gamma(s_j)-\Gamma(s_i)\bigr) \\
\ang\bigl(\Gamma(s_j)-\Gamma(s_i),\Gamma(s_j)-x_j\bigr)
\qquad\mbox{and}\qquad \ang\bigl(\Gamma(s_j)-x_j,v_j\bigr)
\end{gather*}
are smaller than $1/8$, so that the interiors of $K_i$ and $K_j$
certainly are separated by (say) the plane passing through the
midpoint of $[x_{i+1},x_j]$ and perpendicular to
$\Gamma(s_j)-\Gamma(s_i)$, i.e. to the axis of~$K$. This proves
Claim~2.

\bigskip

}
\smallskip

We are now in the position to define the ambient isotopy between
$\gamma$ and $P_\gamma $. Note that $F\colon S_L \to \R^3$ given
by
\[
F(t):= [x_i,x_{i+1}] \cap D_i(\Gamma(t)) \qquad\mbox{for \, $t\in
[t_i,t_{i+1})$, $i=1,\ldots, N$}
\]
is a well defined homeomorphism, parametrizing $P_\gamma$. The
desired isotopy
\[
H\colon \R^3 \times [0,1] \to \R^3
\]
is equal to the identity on $\R^3\setminus \bigcup_{i=1}^N K_i$,
and on each `double cone' $K_i$ it maps each two-dimensional slice
$D_i(z)$, $z\in K_i$, homeomorphically to itself, keeping the
boundary circle of $D_i(z)$ fixed and moving the point $\Gamma(s)$
along a straight segment on $D_i(\Gamma(s))$ until it hits
$[x_i,x_{i+1}]$. \hfill $\Box$

\medskip\noindent\textbf{Proof of Theorem~\ref{thm:1.2}.}
Abbreviate the maximal energy value
$E:=\max\{\E_q(\Gamma_1),\E_q(\Gamma_2)\}$ of the two simple
arclength parametrized curves $\Gamma_i:S_{L_i}\to\R^3$ of
respective (and a priori possibly quite different) lengths $L_i$,
$i=1,2.$ Recall the assumption that the two curves are close in
Hausdorff-distance:
$\dist_H(\Gamma_1,\Gamma_2)<\delta(q)E^{-1(q-2)}.$

Fix $N=N(q,E)$ so that $L_1/N=:\eta< \frac 13 \delta_2(q)
E^{-1/(q-2)}$, set $\eps:=\eta/50$ and let $t_i:=(i-1)\eta\in
S_{L_1}$ for $i=1,\ldots, N$, and $t_{N+1}:=t_1$. By
Proposition~\ref{poly-1}, $\gamma_1$ is ambient isotopic to the
polygonal line
\[
P_{\gamma_1}:=\sum_{i=1}^N [x_i,x_{i+1}]
\]
where $x_i:=\Gamma_1(t_i)$. Now, for $i=1,\ldots,N$ we set
$w_i:=\Gamma_1'(t_i)$,
$\alpha_i:=\Gamma_1\bigl([t_i,t_{i+1}]\bigr)\subset\gamma_1$, and
introduce the half-spaces $H^+_i:=H^+(x_i,x_i+w_i)$ and
$H_i^-:=\R^3\setminus H_i^+$, which are bounded by affine planes
$P_i:=x_i+w_i^\perp$. Consider the tubular regions
\[
T_i:= H_i^+ \cap H_{i+1}^- \cap B_\eps(\alpha_i) .
\]
Their union contains $\gamma_1=\bigcup\alpha_i$; we clearly have
$T_i\cap T_{i+1}=\emptyset$ as $\alpha_{i+1}\subset H_{i+1}^+$.
Moreover, $T_i\cap T_j=\emptyset$ also when $|i-j|>1$. To see
this, we will use Proposition~\ref{prop:5.1} to prove
\begin{equation}
\label{rigid} \inf\{|\Gamma_1(\tau)-\Gamma_1(\sigma)|\ \colon\
(\sigma,\tau)\in S_{L_1}\times S_{L_1},  \ |\sigma-\tau|\ge \eta\}
\ge \frac 34 \eta.
\end{equation}
Before doing so, let us conclude from \eqref{rigid}: If there
existed a point $z\in T_i\cap T_j$ with $|i-j|>1$, we could find
$\sigma\in [t_i,t_{i+1})$ and $\tau \in [t_j,t_{j+1})$ such that
$|\Gamma(\sigma)-\Gamma(\tau)|\le 2 \eps=\eta/25$ by triangle
inequality, a contradiction to \eqref{rigid}.

To verify \eqref{rigid}, we repeat the trick that has already
been used in the proof of Proposition\ref{prop:strong}.
Notice that \eqref{bilip} applied to $\Gamma_1$ implies
\begin{equation}\label{plus}
|\Gamma_1(\tau)-\Gamma_1(\sigma)|\ge\frac 34 |\tau-\sigma|\ge
\frac 34 \eta \quad\Foa\eta\le |\tau-\sigma|\le 3\eta,
\end{equation}
so that the continuously differentiable function $g:S_{L_1}\times
S_{L_1}\to\R$ given by $g(s,t):=|\Gamma_1(s)-\Gamma_1(t)|^2$
attains a positive minimum $g_0>0$ on the compact set
$K_{3\eta}$, where we set $K_\rho:=S_{L_1}\times
S_{L_1}\setminus\{|s-t|<\rho\}$, i.e., there is a pair
of parameters $(s^*,t^*)\in K_{3\eta}$ such that
$g(s,t)\ge g(s^*,t^*)=g_0$ for all $(s,t)\in K_{3\eta}.$
If $|s^*-t^*|=3\eta$ we can apply \eqref{plus} to find
$$
|\Gamma_1(\tau)-\Gamma_1(\sigma)|=\sqrt{g(\tau,\sigma)}
\ge\sqrt{g(s^*,t^*)}=|\Gamma_1(s^*)-\Gamma_1(t^*)|\overset{\eqref{plus}}{\ge}
\frac 34 \eta\quad\Foa (\tau,\sigma)\in K_{3\eta}.
$$
If, on the other hand, $|s^*-t^*|>3\eta$ then by minimality
$\nabla g(s^*,t^*)=0$, which implies that both tangents $\Gamma_1'(s^*)$
and $\Gamma_1'(t^*)$ are perpendicular to the segment $\Gamma_1(s^*)-
\Gamma_1(t^*).$ Thus the intersection
$$
\Gamma_1(S_{L_1})\cap B(\Gamma_1(s^*),2\sqrt{g_0})
\cap B(\Gamma_1(t^*),2\sqrt{g_0})
$$
cannot be contained in the intersection
$C_\phi(\Gamma_1(s^*),\Gamma_1(t^*))\cap
C_\phi(\Gamma_1(t^*),\Gamma_1(s^*)),
$
which according to \eqref{2cones} means that
$$|\Gamma_1(s^*)-\Gamma_1(t^*)|\ge\delta_2(q)E^{\frac{-1}{q-2}}>3\eta,
$$
establishing \eqref{rigid} also in this case.

Assume now that $\dist_{H}(\gamma_1,\gamma_2)< \eps$. We shall
prove that $\gamma_2$ is ambient isotopic to $\gamma_1$; by the
choice of $\eps$, this will mean that Theorem~\ref{thm:1.2}
holds with $\delta_3(q)=\delta_2(q)/150$.

\noindent{\bf Claim.}\, {\it For each $i=1,\ldots,N$ there is a
point
$$
y_i\in P_i\cap\gamma_2\cap B(x_i,2\eps).
$$
}
Without loss of generality we can assume that the curve $\Gamma_1$
is oriented in such a way that
\begin{equation}\label{orientation}
\ang (\Gamma_1'(t_i),v_i)<\frac 18 \quad\AND\quad\ang(\Gamma_1'(t_i),v_{i-1})<\frac 18\quad\Foa i=1,\ldots,N,
\end{equation}
that is, each tangent $\Gamma_1'(t_i)$ points into the set
$K_i:=K(x_i,x_{i+1})=K(\Gamma_1(t_i),\Gamma_1(t_{i+1})),$
which readily implies for the hyperplanes $P_i\perp\Gamma_1'(t_i)$,
$i=1,\ldots,N$,
$$
\ang (P_i,v_i)\ge\ang (P_i,\Gamma_1'(t_i))-\ang(\Gamma_1'(t_i),v_i)
>\frac{\pi}{2}-\frac 18,
$$
and similarly $\ang(P_i,v_{i-1}) > \frac{\pi}{2}-\frac 18$.

Indeed, according to \eqref{2cones}
$$
\Big[\gamma_1\cap B(x_i,2|v_i|) \cap B(x_{i+1},2|v_i|)\cap
H^+(x_i,x_{i+1})\cap H^+(x_{i+1},x_i)\Big] \subset K_i,
$$
which implies that the tangent direction of the curve $\Gamma_1$
at $x_i$ cannot deviate too much from the straight line through
$x_i$ and $x_{i+1}$. The inequalities in \eqref{orientation}
provide a quantified version of this fact.

Since $\dist_H(\gamma_1,\gamma_2)<\eps$ we find three points
$$
z_i\in\gamma_2\cap B(x_i,\eps),\quad
z_{i+1}\in\gamma_2\cap B(x_{i+1},\eps)\quad\AND\quad
z_{i-1}\in\gamma_2\cap B(x_{i-1},\eps)\quad\Foa i=1,\ldots,N.
$$
If $z_i\in P_i$ we set $y_i:=z_i$, and we are done. Else we know
that $z_i\in H^+_i\setminus P_i$ or that $z_i\in H_i^-.$
In the first case we will work with the two points $z_i$ and $z_{i-1}$,
in the second with $z_i$ and $z_{i+1}$ in the same way, so let us
assume the second situation $z_i\in H_i^-.$ We know that $z_{i+1}
\in H^+\setminus P_i$ since by \eqref{bilip}
$$
\dist(z_{i+1},H^-_i)\ge\dist (x_{i+1},H_i^-)-\eps\overset{\eqref{bilip}}{\ge}
\left(\frac 34 - \frac 1{50}\right)\eta >0.
$$
On the other hand, $z_i$ and $z_{i+1}$ are not too far apart,
$$
\rho_i:=|z_i-z_{i+1}|\le|z_i-x_i|+|x_i-x_{i+1}|+|x_{i+1}-z_{i+1}|<
2\eps+\eta<\delta_2(q)E^{-\frac 1{q-2}},
$$
so that we can infer from \eqref{2cones} applied to the points
$x:=z_i$ and $y:=z_{i+1}$ that
\begin{equation}\label{PLUS}
\gamma_2\cap B(z_i,2\rho_i)
\cap B(z_{i+1},2\rho_i)\cap H^+(z_i,z_{i+1})\cap H^+(z_{i+1},z_i)
\subset K(z_i,z_{i+1}).
\end{equation}
We will show that
\begin{equation}\label{PLUSPLUS}
\Big[ K(z_i,z_{i+1})\cap P_i\Big] \subset B(x_i,2\eps).
\end{equation}
Notice that $K(z_i,z_{i+1})\setminus P_i$ consists of two
components, one containing $z_i\in\gamma_2$, and the other
one containing $z_{i+1}\in\gamma_2,$ which implies
that the intersection in \eqref{PLUSPLUS} is not empty.
Since $\gamma_2$ connects $z_i$ and $z_{i+1}$
 by \eqref{PLUS} within the set $K(z_i,z_{i+1})$,
the inclusion in \eqref{PLUSPLUS} yields the desired curve point
$$
y_i\in P_i\cap\gamma_2\cap B(x_i,2\eps)\quad\Foa i=1,\ldots,N,
$$
thus proving the claim.

\medskip

To prove \eqref{PLUSPLUS} we first estimate the angle $\ang(z_{i+1}-
z_i,v_i)$ by the largest possible angle between a line tangent
to both  $B(x_i,\eps)$ and $B(x_{i+1},\eps)$ and the line connecting
the centers $x_i,$ $x_{i+1}$:
$$
\ang(z_{i+1}-z_i,v_i)\le \arcsin\frac{\eps}{|v_i|/2},
$$
so that
\begin{eqnarray*}
\ang(z_{i+1}-z_i,\Gamma_1'(t_i)) & < & \frac 18
+\arcsin\frac{2\eps}{|v_i|} \\
& \overset{\eqref{bilip}}{<}& \frac 18
+\arcsin\frac{2\eta/50}{3\eta/4}  <  \frac 15.
\end{eqnarray*}
%whence $\ang (P_i,z_{i+1}-z_i) > \frac{\pi}{2}-\frac 15.$
Now, let $\tilde{z_i}$ be the orthogonal projection of $z_i$ onto
$P_i$. Since $\ang (\tilde{z}_i-z_i,z_{i+1}-z_i) = \ang
(\Gamma'(t_i),z_{i+1}-z_i) < \frac 15$, it is easy to see that $
K(z_i,z_{i+1}) \cap P_i \subset B(\tilde{z}_i,\tilde h) \cap P_i $
where
$$
\tilde{h}\le|z_i-\tilde{z}_i| \tan\Big(\frac 15 + \frac 18\Big)\le
\eps\tan\Big(\frac{8+5}{40}\Big)<\frac {\eps}2
$$
(see Figure 1 below), which establishes $K(z_i,z_{i+1})\cap
P_i\subset B(x_i,2\eps)$ and hence \eqref{PLUSPLUS}.

\begin{wrapfigure}[19]{l}[0cm]{6.5cm}
\vspace{2cm}
\includegraphics*[totalheight=4.2cm]{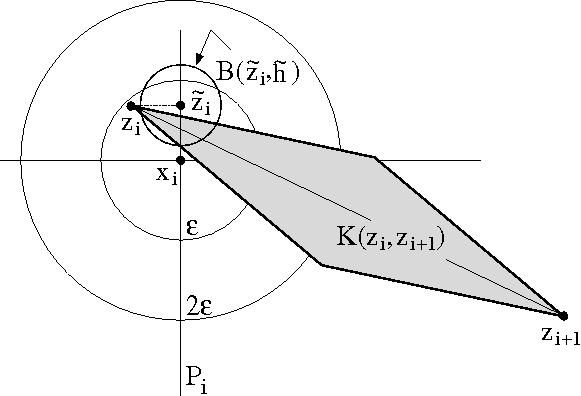}
\vspace{1cm}

{ \footnotesize\baselineskip 11pt {\bf Figure 1.} The intersection
of the doubly conical region $K(z_i,z_{i+1})$ with the plane $P_i$
is contained in the ball $B(\tilde{z}_i,\tilde{h})\subset
B(x_i,2\epsilon)$.

}
\end{wrapfigure}

%For the proof of \eqref{PLUSPLUS} it suffices to replace the point
%$z_i$ by one of the vertices $\tilde{z}_i$ of the cube with edge
%length $2\eps$ centered at $x_i$ (circumscribing the ball
%$B(x_i,\eps)$), and analyzing the worst case situation where a
%doubly conical region $\tilde{K}_i$ congruent to $K(z_i,z_{i+1})$,
%with $\tilde{z}_i$ as one tip point and some point
%$\tilde{z}_{i+1}\in H^+\setminus P_i$ as the other tip, and with
%maximal tilt angle
%$\ang(\tilde{z}_{i+1}-\tilde{z}_i,\Gamma_1'(t_i)) =\frac 15$; see
%Figure ????????????. The height $\tilde{h}$ of the intersection
%$\tilde{K}_i\cap P_i$ can be estimated as

\heikodetail{

\bigskip

Pavel's original argument:

Each ball $B(\Gamma_1(s_i),\eps)$, where $s_i=(t_i+t_{i+1})/2$
contains a point $z_i\in \gamma_2$. We have $|z_i-z_{i+1}|<\eta +
2\eps <2\eta < \delta_2(q) E^{-1/(q-2)}$ so that
Proposition~\ref{prop:5.1} applied to $\gamma_2$ allows us to
conclude that for each $i=1,\ldots,N$ there is a point
\[
y_i=\Gamma_2(\tau_i)\ \in\ B(x_i,\eps) \cap P_i \cap \gamma_2\, .
\]

\bigskip

}

\medskip

Since $|y_i-y_{i+1}|< \eta + 4\eps< 3\eta <\delta_2(q)
E^{-1/(q-2)}$, the curve $\gamma_2$ is ambient isotopic to the
polygonal curve $P_{\gamma_2} =\bigcup_{i=1}^N [y_i,y_{i+1}]$. To
finish the proof of Theorem~\ref{thm:1.2}, it is now sufficient to
check that  $P_{\gamma_1} $ and $P_{\gamma_2} $ are
combinatorially equivalent.

Since the sets $T_i$ are pairwise disjoint, we have
\[
\mathrm{conv} (x_i,x_{i+1},y_i,y_{i+1}) \cap P_{\gamma_1} =
[x_i,x_{i+1}].
\]
This guarantees that all steps in the construction that follows
involve legitimate $\Delta$ and $\Delta^{-1}$-moves. The first
step, taking place in $\overline{T}_1$, is to replace $[x_1,x_2]$
by the union of $[x_1,y_1]$ and $[y_1,x_2]$, and then to replace
$[y_1,x_2]$ by the union of $[y_1,y_2]$ and $[y_2,x_2]$. Next we
perform one $\Delta^{-1}$ and one $\Delta$-move in each of the
$\overline{T}_j$ for $j=2,\ldots, N-1$, replacing first
$[y_j,x_j]$ and $[x_{j},x_{j+1}]$ by $[y_j,x_{j+1}]$, and next
trading $[y_j,x_{j+1}]$ for the union of  $[y_j,y_{j+1}]$ and
$[y_{j+1},x_{j+1}]$. Finally, for $j=N$ we perform two
$\Delta^{-1}$-moves: first replace $[y_N,x_N]$ and $[x_{N},x_{1}]$
by $[y_N,x_1]$, and then replace $[y_N, x_1]$ and $[x_1,y_1]$
(which has been added at the very beginning of the construction)
by $[y_N,y_1]$. This concludes the whole proof.\qed

\section{Bootstrap: optimal regularity of $\Gamma'$}

\label{last}

\setnumbers

In this section, we show how to derive Theorem~\ref{thm:1.3}.
The overall idea is similar to the one in \cite[Section
6]{ssvdm-triple} but here the proof is a little bit less involved.

Assume that $\Gamma$ is 1--1, $\Gamma'\in C^\kappa$, $\kappa =
(q-2)/(q+4)$. Restricting $\Gamma$ to a sufficiently short
interval $I$ in $[0,L]$, and rotating the coordinate system if
necessary, we may assume that the first component $\Gamma_1'$ of
the tangent vector satisfies $\Gamma_1'\ge 0.99$ on $I$ and
$|\Gamma_i'|\approx 0$ on $I$ for all $i=2,\ldots,n$. In fact, to
achieve such control of $\Gamma'$ on $I$ it is enough to assume
that
\[
|I|\le \delta_4(q) \E_q(\Gamma)^{-1/(q-2)}
\]
for some $\delta_4(q)>0$ sufficiently small; the desired
control of $\Gamma'$ follows then from Proposition~\ref{prop:5.1}.

Let
\begin{equation}
\label{max-osc} \Phi(t):= \sup_{{\scriptstyle J\subset
I}\atop{\scriptstyle \mathscr{L}^1(J)\le t}} \left( \osc_J
\Gamma'\right) \qquad\mbox{for $|t|\le \mathscr{L}^1(I) $}
\end{equation}
(here, $J$ denotes an arbitrary subinterval of $I$). We shall show
that for every $u,v\in I$, $u<v$,
\begin{equation}
\label{c-lambda} |\Gamma'(u)-\Gamma'(v)| \le 2
\Phi\Bigl(\frac{|u-v|}N\Bigr)
 + 100
  K_0|u-v|^\lambda,
\end{equation}
where $\lambda=1-2/q$, $N=N(q)>8$ is a large number such that
$2/N^\kappa<1/2$, and
\[
K_0 := \left(N^2 \int_u^v\int_u^v r^{-q} \, ds\, dt\right)^{1/q}\,
.
\]
Once \eqref{c-lambda} is established, we can iterate it to get rid
of the first term on the right hand side of \eqref{c-lambda} and
prove that
\begin{equation}
\label{point}
 |\Gamma'(u)-\Gamma'(v)| \le c_3(q)
\left(\int_u^v\int_u^v r^{-q} \, ds\, dt\, \right)^{1/q}
|u-v|^{1-2/q}\, .
\end{equation}
The argument that shows that \eqref{c-lambda} yields
\eqref{point} is technical but relatively easy; similar
reason\-ings are well known in the theory of PDE (e.g. when one
deals with various Campanato--Morrey estimates). Similar arguments
are described in more detail in our papers
\cite[Section~6]{ssvdm-triple} (see the Remark that follows the
statement of Lemma~6.1 there) and \cite[Section 6]{mengersurf}.
The reader is invited to fill in the computational details or to
consult \cite{ssvdm-triple,mengersurf}.

\bigskip\noindent\textbf{Proof of \eqref{c-lambda}.} We fix
$u<v\in I$ and set
\begin{eqnarray*}
Y_0 & :=& \{ s\in [u,v] \, \colon \H^1(Y_1(s))\ge 2|u-v|/N \}\, ,\\
Y_1(s) & := & \{ t\in [u,v] \, \colon 1/r(\Gamma(s),\Gamma(t))\ge
K_0 |u-v|^{-2/q} \}\, .
\end{eqnarray*}

The reader should think of the parameters in $Y_0$ and $Y_1(s)$ as
`bad' ones. Here is a word of informal explanation. Suppose that a
curve  is just $C^{1,\lambda}$ for $\lambda=1-2/q$ and not
smoother, say like the graph of $x\mapsto |x|^{2-2/q}$ near zero.
We would then expect that a typical point $\Gamma(v)$ can be
roughly at the distance $d^{1+\lambda}$ from the tangent line at
$\Gamma(u)$ when $|\Gamma(u)-\Gamma(v)|\approx d$ or, equivalently
for a flat graph over some interval, $|u-v|\approx d$. But then
$1/r$ at these two points would not exceed a constant multiple of
$d^{1+\lambda}/d^2\approx |u-v|^{-2/q}$  by the explicit formula
\eqref{1.1} for the radius $r$. As we know nothing about the
existence of $\Gamma''$, there are no a priori upper bounds for
$1/r$ that we might use. However, it is illustrative to look at
the sets of points where the model bound $1/r\lesssim
|u-v|^{-2/q}$ is violated. It will turn out that there are `not
too many' such points at all scales, and this will be enough to
conclude.

Set also
\[
E(u,v) := \int_u^v\int_u^v r^{-q} \, ds\, dt\, .
\]
We have
\begin{eqnarray*}
E(u,v) & \ge & \int_{Y_0}\int_{Y_1(s)} r^{-q} \, dt\, ds \\
& \ge & \H^1(Y_0) \cdot \frac{2|u-v|}{N} \, \cdot\, K_0^q
|u-v|^{-2}
 \ = \ \H^1(Y_0) \cdot \frac{2N}{|u-v|}\cdot
 E(u,v)\ ,
\end{eqnarray*}
so that
\[
\H^1(Y_0) \le \frac{|u-v|}{2N}\, .
\]
Now, select $s\in [u,v]\setminus Y_0$ and $t\in [u,v]\setminus
Y_0$ such that
\[
\max(|u-s|,|t-v|) < \frac{|u-v|}N.
\]
By the triangle inequality,
\begin{eqnarray*}
|\Gamma'(u)-\Gamma'(v)|& \le &
|\Gamma'(u)-\Gamma'(s)|+|\Gamma'(s)-\Gamma'(t)|+|\Gamma'(t)-\Gamma'(v)|\\
& \le & 2 \Phi\Bigl(\frac{|u-v|}N\Bigr) +
|\Gamma'(s)-\Gamma'(t)|\, .
\end{eqnarray*}

If the tangent lines $\ell(s)$ and $\ell(t)$ are parallel, we have
$\Gamma'(s)=\Gamma'(t)$ and there is nothing more to prove. Thus,
let us assume that $\ell(s)$ and $\ell(t)$ are not parallel and
proceed to estimate $|\Gamma'(s)-\Gamma'(t)|$.

\medskip

Let $G:= [u,v] \setminus \bigl(Y_1(s)\cup Y_1(t)\bigr)
$. By
definition of $Y_1(\cdot)$ and choice of $s,t$, we have
\begin{equation}
\label{meas-G} \H^1(G)>  |u-v|\left(1-\frac 4N\right) >
\frac{|u-v|}{2}\, .
\end{equation}
If $\sigma\in G$, then by definition of $Y_1(s)$ and of the
tangent-point radius (see \eqref{1.1}) we obtain
\begin{equation}\label{h0} \dist
(\Gamma(\sigma), \ell(s)) <  \frac 12 K_0 |u-v|^{-2/q}
|\Gamma(\sigma)-\Gamma(s)|^2 \le \frac 12 K_0 |u-v|^{2-2/q}=:h_0\,
.
\end{equation}
A similar inequality is satisfied by the distance of
$\Gamma(\sigma)$ to the other line, $\ell(t)$.

\medskip
Now, let $H=\mathrm{span}\, (\Gamma'(s),\Gamma'(t)) \subset
\bbbr^n$ be the two-dimensional plane spanned by the two tangent
vectors $\Gamma'(s)$ and $\Gamma'(t)$. Choose two points $p_1\in
\ell(s)$ and $p_2\in \ell(t)$ such that $|p_1-p_2|= \dist
(\ell(s),\ell(t))$ and let $x:=(p_1+p_2)/2$. (If $\ell(s)$ and
$\ell(t)$ intersect, $x=p_1=p_2$ is their common point; otherwise,
the segment $J(s,t):= [p_1,p_2]$ is perpendicular to each of these
 two lines and $x$ is its midpoint.)

Let $P=x+H$. Then $\dist (\ell(s),P) = \dist(\ell(t),P) =
|p_1-p_2|/2$. Let $\pi_P$ be the orthogonal projection onto $P$
and let
\[
l_1 := \pi_P(\ell(s))\, ,\qquad l_2 := \pi_P(\ell(t))\, .
\]
The lines $l_1,l_2$ intersect at $x\in P$. Note that since $G$ is
nonempty by \eqref{meas-G}, we have in fact by virtue of
\eqref{h0}
\[
|p_1-p_2|=2|x-p_1|\le 2 h_0,
\]
and
\[
\dist(\Gamma(\sigma), l_i) \le 2h_0, \qquad i=1,2, \quad \sigma\in
G\, .
\]
\heikodetail{

\bigskip

$$
\dist(\ell(s),\ell(t))\le\dist (\Gamma(\sigma),\ell(s))+
\dist (\Gamma(\sigma),\ell(t))\overset{\eqref{h0}}{\le} 2h_0,
$$
and
$$
\dist(\Gamma(\sigma),l_1)\le\dist(\Gamma(\sigma),\ell(s))+
\dist(l_1,\ell(s))\le
\dist(\Gamma(\sigma),\ell(s))+\dist(\ell(s),P)\le h_0+\frac{|p_1-p_2|}{2}.
$$

\bigskip

}
Thus,
\begin{equation}
Z:=\Gamma(G) \ \subset \ U_{3 h_0}(l_1) \cap U_{3 h_0}(l_2)\, .
\label{tubes}
\end{equation}
Therefore, the projection $\pi_P(Z)$ of $Z$ onto $P$ is contained
in a rhombus $R$ in $P$. The center of symmetry of $R$ is at $x$;
the sides of $R$ are parallel to $l_1$ and $l_2$; its height
equals $6 h_0$ and its acute angle
$$
\gamma_0:= \ang(l_1,l_2) = \ang(\Gamma'(s),\Gamma'(t))
$$
(since $\Gamma_1'\ge 0.99$ on $I$, the angle
$\ang(\Gamma'(s),\Gamma'(t))$ \emph{is\/} acute). The longer
half-diagonal $D$ of $R$ is given by
\begin{equation}
\label{D} D = \frac{6  h_0}{\sin (\gamma_0/2)}\, ,
\end{equation}
and
\[
\pi_P(Z) \, \subset \, R \, \subset B_D(x) \cap P.
\]
Since $D\ge 6 h_0$, invoking \eqref{tubes} and the triangle
inequality we conclude that
$$
Z=\Gamma(G)\, \subset\, B_{2D}(x)\, .
$$
Now, recall that $\Gamma_1'\ge 0.99$ on $I$. Let $t_2=\sup G$
and $t_1=\inf G$. We then have
\begin{eqnarray*}
4D & = &\diam B_{2D}(x) \ge |\Gamma(t_2)-\Gamma(t_1)|\ge
\Gamma_1(t_2)-\Gamma_1(t_1)\\
& = & \int_{t_1}^{t_2}\Gamma'_1(\sigma)\, d\sigma \\
& \ge & 0.99 \H^1(G)\, .
\end{eqnarray*}
Thus
\begin{eqnarray}
\H^1(G) & < & 5D
 =  \frac{30 h_0}{\sin (\gamma_0/2)} \qquad\mbox{by \eqref{D}} \nonumber \\
& \le & \frac{30 \pi h_0}{\gamma_0} \label{G-above}\
\end{eqnarray}
Combining two estimates of $\H^1(G)$, \eqref{meas-G} and
\eqref{G-above}, we obtain
\[
\ang (\Gamma'(s),\Gamma'(t))=\gamma_0\le \frac{60 \pi
h_0}{|u-v|}\overset{\eqref{h0}}{\equiv} 30\pi K_0
|u-v|^{1-2/q} < 100 K_0 |u-v|^{1-2/q}.
\]
This yields the desired estimate of $|\Gamma'(t)-\Gamma'(s)|$. The
proof of the second part of Theorem~\ref{thm:1.3} is now complete.
\hfill $\Box$

\medskip\noindent\textbf{Remark.} To see that the exponent $1-2/q$
is indeed optimal and cannot be replaced by any larger exponent,
we follow the idea given by M. Szuma\'{n}ska in her PhD thesis
\cite{marta}. One has to fix an arbitrary $a\in (2-2/q,2]$ and
consider $\gamma$ that is the graph of $f(x)=x^a$ say on $[0,1]$.
It is possible to check that $\E_q(\gamma)$ is finite; however,
the derivative of the arclength parametri\-zation of $\gamma$ is
not H\"{o}l\-der continuous with any exponent larger than $\beta=a-1$.
Since $\beta$ can be an arbitrary number in $(1-2/q,1]$, the
exponent $1-2/q$ is indeed optimal. We do not give here the
computational details which are somewhat tedious but routine; one
just has to pass from the graph description of $\gamma$ to the
arclength parametrization and use Taylor's formula in estimates.
The key point is that $f'(x)=ax^{a-1}$ is not H\"{o}l\-der continuous
with any exponent larger than $\beta=a-1$, due to its
behaviour near to 0.

%%%%%%%%%%%%%%   fuer Navigationsleiste  %%%%%%%%%5
\pdfbookmark[0]{References}{References}
%%%%%%%%%%%%%%%%%%%%%%%%%%%%%%%%%%%

%%%%%%%%%%%%%%%%%%%%%%%%%%%%%%%%%%%%%%%%%%%%%
%%%%%%%%%%%%%%%   addresses  %%%%%%%%%%%%%%%%
%%%%%%%%%%%%%%%%%%%%%%%%%%%%%%%%%%%%%%%%%%%%%%

\small
\vspace{1cm}
\begin{minipage}{56mm}
{\sc Pawe\l{} Strzelecki}\\
Instytut Matematyki\\
Uniwersytet Warszawski\\
ul. Banacha 2\\
PL-02-097 Warsaw \\
POLAND\\
E-mail: {\tt pawelst@mimuw.edu.pl}
\end{minipage}
\hfill
\begin{minipage}{56mm}
{\sc Heiko von der Mosel}\\
Institut f\"ur Mathematik\\
RWTH Aachen\\
Templergraben 55\\
D-52062 Aachen\\
GERMANY\\
Email: {\tt heiko@}\\{\tt instmath.rwth-aachen.de}
\end{minipage}

\end{document}